\documentclass[11pt, twoside, a4paper]{article}
\usepackage[english]{babel}
\usepackage{amsmath,amssymb,amsthm,epsfig}
\newtheorem{1}{Theorem}
\newtheorem{2}[1]{Theorem}
\newtheorem{3}[1]{Lemma}
\newtheorem{4}[1]{Proposition}
\newtheorem{5}[1]{Corollary}

\newtheorem{definition}{Definition}

\begin{document}

\title{On the Aubry-Mather theory \\ for symbolic dynamics}
\author{E. Garibaldi$^1$ and A. O. Lopes$^2$}
\date{\today}
\maketitle

\begin{abstract}

We propose a new model of ergodic optimization for expanding
dynamical systems: the holonomic setting. In fact, we introduce
an extension of the standard model used in this theory. The formulation
we consider here is quite natural if one wants a meaning
for possible variations of a real trajectory under the forward shift.
In another contexts (for twist maps, for instance), this property appears
in a crucial way.

A version of the Aubry-Mather theory for symbolic dynamics is
introduced. We are mainly interested here in problems related to
the properties of maximizing probabilities for the two-sided
shift. Under the transitive hypothesis, we show the existence of
sub-actions for H\"older potentials also in the holonomic setting.
We analyze then connections between calibrated sub-actions and the
Ma\~n\'e potential. A representation formula for calibrated
sub-actions is presented, which drives us naturally to a
classification theorem for these sub-actions. We also investigate
properties of the support of maximizing probabilities.

\end{abstract}

\vspace{3cm}

{\footnotesize

\noindent $^1$Institut de Math\'ematiques, Universit\'e Bordeaux 1, F-33405 Talence, France. 

\noindent \emph{Eduardo.Garibaldi@math.u-bordeaux1.fr}

\noindent E. Garibaldi was supported by CAPES PhD scholarship.

\vspace{.5cm}

\noindent $^2$Instituto de Matem\'atica, UFRGS, 91509-900 Porto Alegre, Brasil.

\noindent \emph{alopes@mat.ufrgs.br} 

\noindent A. O. Lopes was partially supported by CNPq, PRONEX -- Sistemas Din\^amicos, Instituto 
do Mil\^enio, and is beneficiary of CAPES financial support.

}

\newpage

\textbf{{\large 1. The Holonomic Condition}}

\vspace {.5cm}

Consider $ X $ a compact metric space. Given a continuous
transformation $ T: X \to X $, we denote by $ \mathcal M_T $ the
convex set of $T$-invariant Borel probability measures. As usual,
we consider on  $ \mathcal M_T $ the weak* topology.

The triple $ (X, T, \mathcal M_T) $ is the standard model used in ergodic optimization. 
Thus, given a potential $ A \in C^0(X) $,
one of the main objectives is the characterization of maximizing
probabilities, that is, the probabilities belonging to
$$\left \{\mu \in \mathcal M_T: \int_X A(x) \; d\mu(x) = \max_{\nu
\in \mathcal M_T} \int_X A(x) \; d\nu(x) \right \}. $$ Several
results were obtained related to this maximizing question, among
them \cite{BLT, Bousch1, Bousch2, CLT, HY, Jenkinson1, Jenkinson2,
LT1}. For maximization with constraints see \cite{Garibaldi, GL,
LT2}. Naturally, if we change the maximizing notion for the
minimizing one, the analogous properties will be true.

Our focus here will be on symbolic dynamics. So let
$ \sigma: \Sigma \to \Sigma $ be a one-sided subshift of finite type
given by a $ r \times r $ transition matrix $ \mathbf M $. More precisely, we have
$$ \Sigma = \left\{ \mathbf x \in \{1, \ldots, r\}^{\mathbb N} : \mathbf M (x_j, x_{j+1}) = 1 \text{ for all } j \geq 0 \right\} $$
and $ \sigma $ is the left shift acting on $ \Sigma $, $ \sigma(x_0, x_1, \ldots) = (x_1, x_2, \ldots) $.
Remind that, fixed $ \lambda \in (0, 1) $, we consider $ \Sigma $ with the
metric $ d (\mathbf x, \bar{\mathbf x}) = \lambda^k $, where
$ \mathbf x = (x_0, x_1, \ldots), \bar{\mathbf x} = (\bar x_0, \bar x_1, \ldots) \in \Sigma $ and
$ k = \min \{j: x_j \ne \bar x_j \} $. 

In this particular situation,
given a continuous potential $ A: \Sigma \to \mathbb{R} $, one should be \emph{a priori} interested in $A$-maximizing
probabilities for the triple $ (\Sigma, \sigma, \mathcal M_\sigma) $.

Nevertheless, this standard model of ergodic optimization has a
main difference to the twist maps theory or to the Lagrangian
Aubry-Mather problem: the dynamics of the shift is not defined
(via a critical path problem) from the potential to be maximized.
In similar terms, in the usual shift standard model, the notion of 
maximizing segment is not present. One would like to have
small variations of a optimal trajectory, by means of a path which
is not a true trajectory, but a small variation of a real
trajectory of the dynamical system. We will describe a model of
ergodic optimization for subshifts of finite type where the
concept of maximizing segment can be introduced: the holonomic
setting. In Aubry-Mather theory for Lagrangian systems (continuous
or discrete time), the set of holonomic probabilities has been
considered before by Ma\~n\'e, Mather, Contreras and Gomes. Main
references on these topics are \cite{Bangert, CI, Fathi, Gomes,
Mane}.

In order to define the holonomic model of ergodic optimization,
we introduce the dual subshift  $ \sigma^*: \Sigma^* \to \Sigma^* $ using as transition
matrix the transposed $ \mathbf M^T $. In clear terms, we consider thus the space 
$$ \Sigma^* = \left\{ \mathbf y \in \{1, \ldots, r\}^{\mathbb N} : \mathbf M (y_{j + 1}, y_j) = 1 \text{ for all } j \geq 0 \right\} $$
and the shift $ \sigma^*(\ldots, y_1, y_0) = (\ldots, y_2, y_1) $.
It is possible, in this way, to identify the space of the dynamics $ (\hat \Sigma, \hat \sigma) $,
the natural extension of $ (\Sigma, \sigma) $, with a subset of
$ \Sigma ^ * \times \Sigma $. In fact, if $ \mathbf y = (\ldots, y_1, y_0) \in \Sigma^* $ and
$ \mathbf x = (x_0, x_1, \ldots) \in \Sigma $, then $ \hat \Sigma $ will be the set of points
$ (\mathbf y, \mathbf x) = (\ldots, y_1, y_0 | x_0, x_1, \ldots) \in \Sigma^* \times \Sigma $ such that
$ (y_0, x_0) $ is an allowed word, namely, such that $ \mathbf M(y_0, x_0) = 1 $.

We define then the transformation $ \tau: \hat \Sigma \to \Sigma $ by
$$ \tau(\mathbf y, \mathbf x) = \tau_{\mathbf y}(\mathbf x) = (y_0, x_0, x_1, \ldots). $$
Note that $\hat \sigma^{-1} (\mathbf y, \mathbf x) = (\sigma^*(\mathbf y), \tau_{\mathbf y}(\mathbf x)) $.

Let $ \mathcal M $ be the convex set of probability measures over the
Borel sigma-algebra of $ \hat \Sigma $.

\begin{definition}
In an analogous way to \cite{Gomes}, we
consider the convex compact subset
$$ \mathcal M_0 = \left \{\hat \mu \in \mathcal M: \int_{\hat \Sigma} f(\tau_{\mathbf y}(\mathbf x)) \; d\hat\mu(\mathbf y,\mathbf x) =
\int_{\hat \Sigma} f(\mathbf x) \; d\hat\mu(\mathbf y, \mathbf x) \; \; \forall \; f \in C^0(\Sigma) \right \}. $$
A probability $ \hat \mu \in \mathcal M_0 $ will be called holonomic.
\end{definition}

Note that $ \mathcal M_{\hat \sigma} \subset \mathcal M_0 $. It is
also not difficult to verify that, whenever $ \mu ^ * \times \mu \in \mathcal M_0 $,
we  have $ \mu \in \mathcal M_\sigma $. Moreover, if
$ \hat \mu \in \mathcal M_0 $, then $ \hat \mu \circ \pi_1^{-1} \in \mathcal M_\sigma $,
where $ \pi_1: \hat \Sigma \to \Sigma $ is the canonical projection. Indeed, if  $ f \in C^0(\Sigma) $,
then
\begin{flushleft}
${\displaystyle \int_{\Sigma} f\circ\sigma(\mathbf x) \; d(\hat\mu\circ\pi_1^{-1}) (\mathbf x) =
\int_{\hat\Sigma} f\circ\sigma(\mathbf x) \; d\hat\mu(\mathbf y, \mathbf x) =} $\\
\vspace {.1cm}
\raggedleft {${\displaystyle = \int_{\hat\Sigma} f\circ\sigma(\tau_{\mathbf y}(\mathbf x)) \; d\hat\mu(\mathbf y, \mathbf x)
= \int_{\hat\Sigma} f(\mathbf x) \; d\hat\mu(\mathbf y, \mathbf x) = \int_{\Sigma} f(\mathbf x) \; d(\hat\mu\circ\pi_1^{-1}) (\mathbf x)} $.}
\end{flushleft}
However, $ \mathcal M_0 $ does not contain just
$\hat\sigma$-invariant probabilities. In fact, if $\mathbf x \in
\Sigma $ is a periodic point of period $ M $, fix any subset
$ \{\mathbf y^0, \ldots, \mathbf y^{M - 1} \} \subset \Sigma^* $
with $ y_0^j = x_{M - 1 + j} $ for $ 0 \leq j \leq M - 1 $.
It is easy to see that
$$ \hat \mu = \frac{1}{M} \sum_{j = 0}^{M - 1} {\delta_{\mathbf y^j} \times \delta_{\sigma^j(\mathbf x)}} \in \mathcal M_0. $$

For the ergodic optimization problem, there is very little difference
(in a purely abstract point of view) in relation to which convex compact set of
probability measures over the Borel sigma-algebra is made the
maximization. In fact, an adaptation of the proposition 10 of
\cite{CLT} assures that, when considering a convex compact subset
$ \mathcal N \subset \mathcal M $, a generic H\"older potential
admits a single maximizing probability in $ \mathcal N $.

Taking a continuous application $ A: \hat \Sigma \to \mathbb R $, a natural situation is then to
formulate the maximization problem over the set $ \mathcal M_0 $.

\begin{definition}
Given a potential $ A \in C^0(\hat \Sigma) $, denote
$$ \beta_A = \max_{\hat \mu \in \mathcal M_0} \int_{\hat \Sigma} A(\mathbf y, \mathbf x) \; d\hat\mu(\mathbf y, \mathbf x). $$
\end{definition}

We point out that sometimes, even if one is interested just in
the problem for a H\"older potential $ A: \Sigma \to \mathbb R $, one has to go to the
dual problem and consider the dual potential $ A^*:  \Sigma^* \to \mathbb R $.
This happens, for instance, when someone is trying to analyze a large deviation
principle for the equilibrium probabilities associated to the family of H\"older potentials
$ \{ t A \}_{t > 0} $ (see \cite{BLT}).

Actually, the maximization problem over $ \mathcal M_{\hat \sigma} $ is not so interesting, 
because any H\"older potential $ A : \hat \Sigma \to \mathbb R $ is cohomologous 
to a potential that depends
just on future coordinates (see, for instance, \cite{PP}). In this
case, the problem can be in principle analyzed in the standard
model, that is, over $ \mathcal M_\sigma $.

Furthermore, in order to analyze maximization of the integral
of a potential $ A \in C^0(\Sigma) $, no new maximal value will be
found, because
$$ \max_{\hat \mu \in \mathcal M_0}  \int_{\hat \Sigma} A(\mathbf x) \; d\hat\mu(\mathbf y, \mathbf x) =
\max_{\mu \in \mathcal M_\sigma} \int_{\Sigma} A(\mathbf x) \; d\mu(\mathbf x). $$
Indeed, the correspondence $ \hat \mu \in \mathcal M_0 \mapsto \hat \mu \circ \pi_1^{-1} \in \mathcal M_\sigma $
preserves the integration on $ C^0(\Sigma) $ and the
same property is verified by the correspondence $ \mu \in \mathcal M_\sigma \mapsto \mu \circ \pi_1 \circ \hat \sigma^{-1} \in \mathcal M_0 $.

Therefore, we could say that the holonomic model of ergodic optimization $ (\hat \Sigma, \hat \sigma, \mathcal M_0) $
is an extension of the standard model $ (\Sigma, \sigma, \mathcal M_\sigma) $.

This paper is part of the first author's PhD thesis \cite{Garibaldi}.
We will be interested here in the maximization question over $ \mathcal M_0 $ and,
if possible, in some properties that one can get for the problem over
$ (\Sigma, \sigma) $. In the  section 2, we will show the dual identity
$$ \beta_A =
\inf_{f \in C^0(\Sigma)} \max_{(\mathbf y, \mathbf x) \in \hat \Sigma} [A(\mathbf y, \mathbf x) + f(\mathbf x) - f(\tau_{\mathbf y}(\mathbf x))]. $$

We will then analyze the problem of finding a function $ u \in C^0(\Sigma) $
which realizes the infimum of the previous expression, that is, a sub-action for $ A $.

\begin{definition}
A sub-action $ u \in C^0(\Sigma) $ for the
potential $ A \in C^0(\hat \Sigma) $ is a function satisfying, for
any $ (\mathbf y, \mathbf x) \in \hat \Sigma $,
$$ u(\mathbf x) \le u(\tau_{\mathbf y}(\mathbf x)) - A(\mathbf y, \mathbf x) + \beta_A. $$
\end{definition}

Assuming the dynamics $ (\Sigma, \sigma) $ is topologically mixing and the potential $ A $ is H\"older, 
we will show in section 3 the existence of a H\"older sub-action of maximal character.
Furthermore, under the transitivity hypothesis, for a potential $\theta$-H\"older,
we will show that we can always find a calibrated sub-action $ u \in C^\theta(\Sigma) $.

\begin{definition}
A calibrated sub-action $ u \in C^0(\Sigma) $
for $ A \in C^0(\hat \Sigma) $ is a function satisfying
$$ u(\mathbf x) = \min_{\mathbf y \in \Sigma_{\mathbf x}^*} [u(\tau_{\mathbf y}(\mathbf x)) - A(\mathbf y, \mathbf x) + \beta_A], $$
where, for each point $ \mathbf x \in \Sigma $, we denote by $ \Sigma_{\mathbf x}^* $ the subset of elements
$ \mathbf y \in \Sigma^* $ such that $(\mathbf y, \mathbf x) \in \hat \Sigma $.
\end{definition}

In the transitive context, we will introduce in section 4 the
Ma\~n\'e potential $ S_A: \Sigma \times \Sigma \to \mathbb R \cup \{+\infty \} $
(the terminology is borrowed from Aubry-Mather theory). Thus, we will establish a
family of H\"older calibrated sub-actions, namely, $ \{S_A(\mathbf x, \cdot) \}_{\mathbf x \in \Omega(A)} $,
where $ \Omega(A) $ denotes the set of non-wandering points with respect to the potential
$ A \in C^\theta(\hat \Sigma) $. All these notions will be precisely defined later.
Besides, these concepts already appear in \cite{CLT} for the forward shift setting.

\begin{definition}
We will denote by
$$ \text{\Large $ \mathit m $}_A =
\left \{\hat \mu \in \mathcal M_0 : \int_{\hat \Sigma} A(\mathbf y, \mathbf x) \; d\hat\mu(\mathbf y, \mathbf x) = \beta_A \right \} $$
the set of the $A$-maximizing holonomic probabilities.
\end{definition}

When we investigate the connections between sub-actions and the
supports of holonomic probabilities, the $A$-maximizing holonomic
probability notion is of great importance. One of the main results
of section 5 is the representation formula for calibrated
sub-actions. More specifically, given a calibrated sub-action $ u $ 
for a potential $ A \in C^\theta(\hat \Sigma) $, the following expression holds
$$ u(\bar{\mathbf x}) = \inf_{\mathbf x \in \Omega(A)} [u(\mathbf x) + S_A(\mathbf x, \bar{\mathbf x})]. $$
Such characterization is analogous to the one obtained for weak KAM solutions in Lagrangian systems (see \cite{Contreras}).
Under the transitivity hypothesis, this representation formula and its reciprocal will describe, by means of an 
isometric bijection, the set of the calibrated sub-actions for a
H\"older potential $ A $. We will show yet that $ \hat \mu \in \text{\Large $\mathit m $}_A $
with $ \hat \mu \circ \pi_1^{-1} $ ergodic implies $ \pi_1(\text{supp}(\hat \mu)) \subset \Omega(A) $. This property will
drive us naturally to other questions like, for instance, the possibility of reducing contact loci.

\vspace {1cm}

\textbf{{\large 2. The Dual Formulation}}

\vspace {.5cm}

We start presenting the main goal of this section.

\begin{1}\label{formuladual}
Given a potential $ A \in C^0(\hat \Sigma) $, we have
$$ \beta_A =
\inf_{f \in C^0(\Sigma)} \max_{(\mathbf y, \mathbf x) \in \hat \Sigma} [A(\mathbf y, \mathbf x) + f(\mathbf x) - f(\tau_{\mathbf y}(\mathbf x))]. $$
\end{1}

One observes that this formula corresponds in Lagrangian Aubry-Mather theory
to the characterization of Ma\~n\'e's critical value (see theorem A of \cite{CIPP}).
Theorem~\ref{formuladual} is just a consequence of the Fenchel-Rockafellar theorem. 
For the standard model $ (X, T, \mathcal M_T) $, a similar 
result was established before (consult, for instance, \cite{CG, Radu}). 
We will present, anyway, the complete proof for the holonomic setting.

First, consider the convex correspondence
$ F: C^0(\hat \Sigma) \to \mathbb R $ defined by $ F(g) = \max (A + g) $.
Consider also the subset
$$ \mathcal C = \{g \in C^0(\hat \Sigma): g(\mathbf y, \mathbf x) =
f(\mathbf x) - f(\tau_{\mathbf y}(\mathbf x)),\,\, \text {for some} \,\,f \in C^0(\Sigma) \}. $$
We establish then a concave correspondence $ G: C^0(\hat \Sigma) \to \mathbb R \cup \{- \infty \} $ taking
$ G(g) = 0 $ if $ g \in \bar{\mathcal C} $ and $ G(g) = - \infty $ otherwise.

Let $ \mathcal S $ be the set of the signed measures over the Borel
sigma-algebra of $ \hat \Sigma $. Remember that the corresponding Fenchel tranforms,
$ F^*: \mathcal S \to \mathbb R \cup \{+ \infty \} $ and $ G^*: \mathcal S \to \mathbb R \cup \{- \infty \} $, are given by
$$ F^* (\hat \mu) = \sup_{g \in C^0(\hat \Sigma)} \left [\int_{\hat \Sigma} g(\mathbf y, \mathbf x) \; d\hat\mu(\mathbf y, \mathbf x) - F(g) \right] \;
\text {and} $$
$$ G^* (\hat \mu) = \inf_{g \in C^0(\hat \Sigma)} \left [\int_{\hat \Sigma} g(\mathbf y, \mathbf x) \; d\hat\mu(\mathbf y, \mathbf x) - G(g) \right]. $$
Denote
$$ \mathcal S_0 = \left \{\hat \mu \in \mathcal S:  \int_{\hat \Sigma} f(\tau_{\mathbf y}(\mathbf x)) \; d\hat\mu(\mathbf y, \mathbf x) =
\int_{\hat \Sigma} f(\mathbf x) \; d\hat\mu(\mathbf y, \mathbf x) \; \; \forall \; f \in C^0(\Sigma) \right \}. $$

\begin{3}\label{transformadas}
Given $ F $ and $ G $ as above, we verify
$$ F^* (\hat \mu) =
\left\{ \begin{array}{ll} {\displaystyle - \int_{\hat \Sigma} A(\mathbf y, \mathbf x) \; d\hat\mu(\mathbf y, \mathbf x)}
& \mbox{if $ \hat \mu \in \mathcal M $} \\ + \infty & \mbox{otherwise} \end{array} \right.
\; \text{ and} $$
$$ G^* (\hat \mu) =
\left\{ \begin{array}{ll} 0 & \mbox{if $ \hat \mu \in \mathcal S_0 $} \\ - \infty & \mbox{otherwise} \end{array} \right.. $$
\end{3}

\begin{proof}
Assume first that $ \hat \mu \in \mathcal S $ is not positive,
that is, $ \hat \mu $ gives a negative value for some Borel set.
Therefore, we can find a sequence of functions $ \{g_j \} \subset C^0(\hat \Sigma, \mathbb R^-) $ such that
$ {\displaystyle \lim \int_{\hat \Sigma} g_j(\mathbf y, \mathbf x) d\hat\mu (\mathbf y, \mathbf x) = + \infty} $. Once $ F(g_j) \leq F(0) < + \infty $,
we have $ F^*(\hat \mu) = + \infty $.

Suppose $\hat \mu \in \mathcal S $ is such that $\hat \mu \ge 0 $
and $\hat \mu (\hat \Sigma) \ne 1 $. In this case, we observe
\begin{eqnarray*}
\sup_{g \in C^0(\hat \Sigma)} \left[ \int_{\hat \Sigma} g(\mathbf y, \mathbf x) \; d\hat\mu(\mathbf y, \mathbf x) - F(g) \right]
& \ge & \sup_{a \in \mathbb R} \left[ \int_{\hat \Sigma} a \; d\hat\mu(\mathbf y, \mathbf x) - F(a) \right] \\
& = & \sup_{a \in \mathbb R} \left[ a(\hat\mu(\hat \Sigma) - 1) - F(0) \right] = +\infty.
\end{eqnarray*}

On the other hand, when we consider $ \hat \mu \in \mathcal M $,
directly from the inequality
$ {\displaystyle \int_{\hat \Sigma} A(\mathbf y, \mathbf x) \; d\hat\mu(\mathbf y, \mathbf x)} +
{\displaystyle \int_{\hat \Sigma} g(\mathbf y, \mathbf x) \; d\hat\mu(\mathbf y, \mathbf x)} \leq F(g) $,
we have
$$ - \int_{\hat \Sigma} A(\mathbf y, \mathbf x) \; d\hat\mu(\mathbf y, \mathbf x) \ge
\sup_{g \in C^0(\hat \Sigma)} \left [\int_{\hat \Sigma} g(\mathbf y, \mathbf x) \; d\hat\mu(\mathbf y, \mathbf x) - F(g) \right]. $$
Once $ F(-A) = 0 $, we get the characterization of $ F^* $.

Now we will consider $ G^* $. If $ \hat \mu \notin \mathcal S_0 $,
there exists a function $ f \in C^0(\Sigma) $ such that
$ {\displaystyle \int_{\hat \Sigma} f(\tau_{\mathbf y}(\mathbf x)) \; d\hat\mu(\mathbf y, \mathbf x) \ne
\int_{\hat \Sigma} f(\mathbf x) \; d\hat\mu(\mathbf y, \mathbf x)} $.
Therefore, we verify
\begin{eqnarray*}
G^*(\hat \mu) & = &
\inf_{g \in \mathcal C} \int_{\hat \Sigma} g(\mathbf y, \mathbf x) \; d\hat\mu(\mathbf y, \mathbf x) \\
& \le &
\inf_{a \in \mathbb R} a \int_{\hat \Sigma} \left [ f(\tau_{\mathbf y}(\mathbf x)) - f(\mathbf x) \right ]\; d\hat\mu(\mathbf y, \mathbf x)
= - \infty.
\end{eqnarray*}

Besides, for $ \hat \mu \in \mathcal S_0 $, clearly $ G^*(\hat \mu) = 0 $.
\end{proof}

Using this lemma, we can show the dual expression of the beta constant $ \beta_A = {\displaystyle \max_{\hat \mu \in \mathcal M_0}
\int_{\hat \Sigma} A(\mathbf y, \mathbf x) \; d\hat\mu(\mathbf y, \mathbf x)} $.

\begin{proof}[Proof of Theorem~\ref{formuladual}]
Once the correspondence $ F $ is Lipschitz, the theorem
of duality of Fenchel-Rockafellar assures
$$ \sup_{g \in C^0(\hat \Sigma)} \left[G(g) - F(g)\right] = \inf_{\hat \mu \in \mathcal S} \left[F^*(\hat \mu) - G^*(\hat \mu)\right]. $$
Thus, by lemma~\ref{transformadas},
$$ \sup_{g \in \mathcal C} \left [- \max_{(\mathbf y, \mathbf x) \in \hat \Sigma} (A + g) (\mathbf y, \mathbf x) \right] =
\inf _{\hat \mu \in \mathcal M_0} \left [- \int_{\hat \Sigma} A(\mathbf y, \mathbf x) \; d\hat\mu(\mathbf y, \mathbf x) \right]. $$
Finally, from the definition of $ \mathcal C $, we get the statement of the theorem.
\end{proof}

Relative maximization is studied in \cite{GL}. In this case,
the dual formula is also true. More specifically, if we introduce a constraint
$ \varphi \in C^0(\hat \Sigma, \mathbb R^n) $ with coordinate functions
$ \varphi_1, \ldots, \varphi_n $, we can then consider an induced map
$ \varphi_* \in C^0(\mathcal M_0, \mathbb R^n) $ given by
$$ \varphi_*(\hat \mu) =
\left (\int_{\hat \Sigma} \varphi_1(\mathbf y, \mathbf x) \; d\hat\mu(\mathbf y, \mathbf x), \ldots,
\int_{\hat \Sigma} \varphi_n(\mathbf y, \mathbf x) \; d\hat\mu(\mathbf y, \mathbf x) \right). $$
Thus, if $ A \in C^0(\hat \Sigma) $, we can immediately define a concave
and continuous function $ \beta_{A, \varphi}: \varphi_*(\mathcal M_0) \to \mathbb R $ by
$$ \beta_{A, \varphi}(h) =
\max_{\hat \mu \in \varphi_*^{-1}(h)} \int_{\hat \Sigma} A(\mathbf y, \mathbf x) \; d\hat\mu(\mathbf y, \mathbf x). $$
Using a little bit more refined argument as \cite{Radu}, we
could demonstrate the dual formula for a beta function
$$ \beta_{A, \varphi}(h) = \inf_{(f,c) \in C^0(\Sigma) \times \mathbb R^n} \max_{(\mathbf y, \mathbf x) \in \hat \Sigma}
(A + f \circ \pi_1 - f \circ \pi_1 \circ \hat \sigma^{-1} - \langle c, \varphi - h \rangle)(\mathbf y, \mathbf x). $$

Nevertheless, the unconstrained dual formula raises a natural question: can we
find functions accomplishing the infimum of the dual expression?
In an equivalent way, is there a function $ u \in C^0(\Sigma) $
such that
$$ A + u \circ \pi_1 - u \circ \pi_1 \circ \hat \sigma^{-1}  \leq \beta_A? $$
As we mentioned at the first section, we call any function $ u $
as above a sub-action for $ A $. This terminology is motivated by
the inequality
$$ A + u \circ \sigma - u \leq \beta_A, $$
which is present at the usual definition of a sub-action $ u $
for the forward shift setting (see \cite{CLT} for
instance). The next sections are mainly dedicated to show the existence of
sub-actions in the holonomic setting.

\vspace {1cm}

\textbf{{\large 3. Sub-actions: Maximality and Calibration}}

\vspace {.5cm}

We start showing not only the existence of sub-actions but, as a matter of fact, the existence of a maximal 
sub-action. To that end, remember that a dynamical system $ (X, T) $ is topologically mixing,
if, for any  pair of non-empty open sets $ D, E \subset X $, there is an integer
$ K > 0 $ such that $ T^k(D) \cap E \neq \emptyset $ for all $ k > K $.

\begin{4}
Consider any topologically mixing
subshift of finite type $ \sigma: \Sigma \to \Sigma $ and a
potential $ A \in C^\theta(\hat \Sigma) $. Then, there exists a 
sub-action  $ u_A \in C^\theta(\Sigma, \mathbb R^-) $ such
that, for any other sub-action $ u \in C^0(\Sigma, \mathbb R^-) $,
we have $ u_A \ge u $.
\end{4}

A sub-action like this one (not necessarily H\"older)
will be called maximal.

\begin{proof}
Without loss of generality, we can assume $ \beta_A = 0 $. Then,
for each $ \mathbf x \in \Sigma $, set
$$ u_A (\mathbf x) = \inf \left \{- \sum_{j = 0}^{k - 1} A(\mathbf y^j, \mathbf x^j):
k \ge 0, \; \mathbf x^0 = \mathbf x, \; \mathbf y^j \in \Sigma_{\mathbf x^j} ^ *, \; \mathbf x^{j + 1} = \tau_{\mathbf y^j}(\mathbf x^j) \right \}. $$
By convention, we assume the sum is zero when $ k = 0 $.

Suppose for a moment that $ u_A $ is a well defined H\"older application. Note that,
if $ \mathbf y^0 = \mathbf y $ and $ \mathbf x^0 = \mathbf x $, then
\begin{eqnarray*}
A(\mathbf y, \mathbf x) & = &
\sum_{j = 0}^k A(\mathbf y^j, \mathbf x^j) - \sum_{j = 0}^{k - 1} A(\mathbf y^{j + 1}, \mathbf x^{j + 1}) \\
& \le & - \sum_{j = 0}^{k - 1} A(\mathbf y^{j + 1}, \mathbf x^{j + 1}) - u_A (\mathbf x).
\end{eqnarray*}
Clearly $ \mathbf x^1 = \tau_{\mathbf y^0}(\mathbf x^0) = \tau_{\mathbf y}(\mathbf x) $.
Thus, since the inequality is true for all $ k \ge 0 $ and any
points $ (\mathbf y^1, \mathbf x^1), \ldots, (\mathbf y^k, \mathbf x^k) \in \hat \Sigma $
such that $\mathbf x^{j + 1} = \tau_{\mathbf y^j}(\mathbf x^j) $, it follows that
$ A(\mathbf y, \mathbf x) \le u_A (\tau_{\mathbf y}(\mathbf x)) - u_A (\mathbf x) $, that is, $ u_A $ is a sub-action for the potential $ A $.

So let us prove that the function $ u_A $ is well defined.
Remember that, when $ \bar{\mathbf x} \in \Sigma $ is a periodic point of period $ k $, if we choose any points
$ \bar{\mathbf y}^j \in \Sigma^* $ satisfying $ \bar y_0^j = \bar x_{k - (j + 1)} $, we obtain
$ {\displaystyle \hat \mu =
\frac{1}{k} \sum_{j = 0}^{k - 1} {\delta_{\bar{\mathbf y}^j} \times \delta_{\sigma^{k - j}(\bar{\mathbf x})}} \in \mathcal M_0} $.
Hence, we immediately verify
$$ - \sum_{j = 0}^{k - 1} A(\bar{\mathbf y}^j, \sigma^{k - j}(\bar{\mathbf x})) =
- k \int_{\hat \Sigma} A(\mathbf y, \mathbf x) \; d\hat\mu(\mathbf y, \mathbf x) \ge 0. $$
Given $ \mathbf x \in \Sigma $, we choose then points $ (\mathbf y^0, \mathbf x^0), \ldots, (\mathbf y^{k - 1}, \mathbf x^{k - 1}) \in \hat \Sigma $
satisfying $ \mathbf x^0 = \mathbf x $ and $ \mathbf x^{j + 1} = \tau_{\mathbf y^j}(\mathbf x^j) $. As $ (\Sigma, \sigma) $ is topologically mixing,
there exists an integer $ K > 0 $ such that, for any $ k > K $, we can find a periodic point $ \bar{\mathbf x} $ of period $ k $ satisfying
$ d(\mathbf x^k, \bar{\mathbf x}) < \lambda^{k - K} $, where $ \mathbf x^k = \tau_{\mathbf y^{k - 1}} (\mathbf x^{k - 1}) $.
Thus, when we put $ \bar{\mathbf y}^j = \mathbf y^j $ for $ K \le j \le k - 1 $, it follows that
$$ \left | \sum_{j = 0}^{k - 1} A(\mathbf y^j, \mathbf x^j) - \sum_{j = 0}^{k - 1} A(\bar{\mathbf y}^j, \sigma^{k - j}(\bar{\mathbf x})) \right |
\le \frac{\text{H\"old}_\theta(A)}{1 - \lambda^\theta} + 2 K \| A \|_0, $$
which assures that $ u_A $ is well defined.

The application $ u_A $ is $\theta$-H\"older. Indeed, fix $ \mathbf x, \bar{\mathbf x} \in \Sigma $ with
$ d(\mathbf x, \bar{\mathbf x}) \leq \lambda $ and consider once more points
$ (\mathbf y^0, \mathbf x^0), \ldots, (\mathbf y^{k - 1}, \mathbf x^{k - 1}) \in \hat \Sigma $ satisfying
$ \mathbf x^0 = \mathbf x $ and $ \mathbf x^{j + 1} = \tau_{\mathbf y^j}(\mathbf x^j) $.
Putting $ \bar{\mathbf x}^0 = \bar{\mathbf x} $ and $ \bar{\mathbf x}^{j + 1} = \tau_{\mathbf y^j}(\bar{\mathbf x}^j) $, we obtain
$$ \left | \sum_{j = 0}^{k - 1} A(\mathbf y^j, \mathbf x^j) - \sum_{j = 0}^{k - 1} A(\mathbf y^j, \bar{\mathbf x}^j) \right |
\le \frac{\text{H\"old}_\theta(A)} {1 - \lambda^\theta} d(\mathbf x, \bar{\mathbf x})^\theta. $$
As the collection of points $ \{ (\mathbf y^j, \mathbf x^j) \} $ was chosen arbitrarily, it follows that
$$ | u_A(\mathbf x) - u_A(\bar{\mathbf x}) | \le \frac{\text{H\"old}_\theta(A)} {1 - \lambda^\theta} d(\mathbf x, \bar{\mathbf x})^\theta. $$

To prove the maximal character of $ u_A $, just observe that, for any
sub-action $ u \in C^0(\Sigma, \mathbb R ^ -) $, we have
$$ u(\mathbf x) \leq u(\tau_{\mathbf y^{k - 1}} (\mathbf x^{k - 1})) - \sum_{j = 0}^{k - 1} A(\mathbf y^j, \mathbf x^j) \le
- \sum_{j = 0}^{k - 1} A(\mathbf y^j, \mathbf x^j) $$
when $ k \ge 0 $, $ \mathbf x^0 = \mathbf x $, $ \mathbf y^j \in \Sigma_{\mathbf x^j}^* $ and $ \mathbf x^{j + 1} = \tau_{\mathbf y^j}(\mathbf x^j) $.
\end{proof}

An interesting question is the existence of a sub-action of minimal character.
Given a potential $ A \in C^\theta(\hat \Sigma) $, a possible approach to this demand
is to introduce the function $ U^{K, \theta}_A \in C^\theta(\Sigma) $ defined by
$$ U^{K, \theta}_A = \inf \{u \in C^\theta(\Sigma): u \; \text{ sub-action for } \; A, \; \text{H\"old}_\theta(u) \le K, \; \max u = 0 \}. $$
The sub-action $ U^{K, \theta}_A $ is in some sense minimal.

In the final section, instead of imposing $ \max u = 0 $, we will 
consider a suitable normalization of sub-actions in order to present a maximal calibrated one. 
We will need however several results before to discuss this special situation. For instance,
the following theorem assures the existence of calibrated sub-actions for any
$\theta$-H\"older potential.

\begin{2}\label{existenciacalibrada}
Let $ \sigma: \Sigma \to \Sigma $ be a transitive subshift of finite type. For each potential
$ A \in C^\theta(\hat \Sigma) $, there exists a function $ u \in C^\theta(\Sigma) $ such that
$$ u(\mathbf x) = \min_{\mathbf y \in \Sigma_{\mathbf x} ^ *} [u(\tau_{\mathbf y}(\mathbf x)) - A(\mathbf y, \mathbf x) + \beta_A]. $$
\end{2}

\begin{proof}
The idea is to obtain a fixed point of a weak contraction as a
limit of fixed points of strong contractions (see \cite{Bousch1, Bousch2}).

Given $ \rho \in (0, 1] $, we define the transformation $ \mathcal L_\rho: C^0(\Sigma) \to C^0(\Sigma) $ by
$$ \mathcal L_\rho(f)(\mathbf x) = \rho \min_{\mathbf y \in \Sigma_{\mathbf x}^*} [f(\tau_{\mathbf y}(\mathbf x)) - A(\mathbf y, \mathbf x)]. $$
Once $ \mathcal L_\rho $ is $\rho$-Lipschitz, consider, when $ 0 < \rho < 1 $, its fixed point $ u_\rho \in C^0(\Sigma) $.

The first fact to be noticed is the equicontinuity of the family
$\{u_\rho \} $. Indeed, note that $ \Sigma^*_{\mathbf x^0} = \Sigma^*_{\bar {\mathbf x}^0} $ when
$ d(\mathbf x^0, \bar {\mathbf x}^0) \leq \lambda $. Hence, if  $ \mathbf y^0 \in \Sigma^*_{\mathbf x^0} $ satisfies
$$ u_\rho (\mathbf x^0) = \rho [u_\rho(\tau_{\mathbf y^0}(\mathbf x^0)) - A(\mathbf y^0, \mathbf x^0)], $$
we obtain
$$ u_\rho (\bar {\mathbf x}^0) \leq \rho [u_\rho(\tau_{\mathbf y^0}(\bar {\mathbf x}^0)) - A(\mathbf y^0, \bar {\mathbf x}^0)]. $$
Therefore, taking $ \mathbf x^1 = \tau_{\mathbf y^0}(\mathbf x^0) $ and
$ \bar {\mathbf x}^1 = \tau_{\mathbf y^0}(\bar {\mathbf x}^0) $, we have the inequality
$$ u_\rho (\bar {\mathbf x}^0) - u_\rho (\mathbf x^0) \leq
\rho [A(\mathbf y^0, \mathbf x^0) - A(\mathbf y^0, \bar {\mathbf x}^0)] + \rho [u_\rho (\bar {\mathbf x}^1) - u_\rho (\mathbf x^1)]. $$
In this way, defining $ \mathbf x^j = \tau_{\mathbf y^{j - 1}} (\mathbf x^{j - 1}) $ and
$ \bar {\mathbf x}^j = \tau_{\mathbf y^{j - 1}} (\bar {\mathbf x}^{j - 1}) $, we continue inductively obtaining
$ \mathbf y^j \in \Sigma^*_{\mathbf x^j} $ such that
$ u_\rho (\mathbf x^j) = \rho [u_\rho(\tau_{\mathbf y^j}(\mathbf x^j)) - A(\mathbf y^j, \mathbf x^j)] $.
As a consequence of this construction, it follows
$$ u_\rho (\bar {\mathbf x}^0) - u_\rho (\mathbf x^0) \leq
\sum_{j = 0}^{k - 1} \rho^{j + 1} [A(\mathbf y^j, \mathbf x^j) - A(\mathbf y^j, \bar {\mathbf x}^j)] + \rho^k [u_\rho (\bar {\mathbf x}^k) - u_\rho (\mathbf x^k)].$$
Thus, we verify
\begin{eqnarray*}
u_\rho (\bar {\mathbf x}^0) - u_\rho (\mathbf x^0)
& \le &
\sum_{j = 0}^{\infty} \rho^{j + 1} [A(\mathbf y^j, \mathbf x^j) - A(\mathbf y^j, \bar {\mathbf x}^j)] \\
& \le &
\text{H\"old}_\theta(A) \sum_{j = 0}^{\infty} \rho^{j + 1} d(\mathbf x^j,\bar {\mathbf x}^j)^\theta \\
& \le &
\text{H\"old}_\theta(A) d(\mathbf x^0,\bar {\mathbf x}^0)^\theta \sum_{j = 0}^{\infty} \rho^{j + 1} \lambda^{j\theta} \\
& = &
\frac{\rho \text{H\"old}_\theta(A)}{1 - \rho\lambda^\theta}d(\mathbf x^0,\bar {\mathbf x}^0)^\theta.
\end{eqnarray*}
We proved that the family $ \{ u_\rho \} $ is uniformly $\theta$-H\"older, in particular it is an equicontinuous family of functions.

The family $ \{ u_\rho \} $ presents also uniformly bounded oscillation. Indeed, given a point $ (\mathbf y, \mathbf x) \in \hat \Sigma $, note that
\begin{eqnarray*}
u_\rho(\mathbf x) - \min u_\rho
& \le &
\rho [u_\rho(\tau_{\mathbf y}(\mathbf x)) - A(\mathbf y, \mathbf x)] - \min \rho [u_\rho \circ \pi_1 \circ \hat \sigma^{-1} - A] \\
& \le &
\rho [\max A - A(\mathbf y, \mathbf x)] + \rho[u_\rho(\tau_{\mathbf y}(\mathbf x)) - \min u_\rho] \\
& \le &
\text{H\"old}_\theta(A) + u_\rho(\tau_{\mathbf y}(\mathbf x)) - \min u_\rho.
\end{eqnarray*}
Since $ (\Sigma, \sigma) $ is transitive, we can define a finite set
$ \{(\mathbf y^j, k_j) \} \subset \Sigma^* \times \mathbb N $ by choosing,
for each pair of symbols $ s, s' \in \{1, \ldots, r\} $, an allowed word
$ (y_{k_j - 1}^j, \ldots, y_0^j) $ such that $ y_{k_j - 1}^j = s' $ and the word $ (y_0^j, s) $ is allowed.
Consequently, given $ \mathbf x \in \Sigma $ with $ x_0 = s $, the inequality
$$ u_\rho(\mathbf x) - \min u_\rho \leq k_j \text{H\"old}_\theta(A) + u_\rho(\tau_{\mathbf y^j}^{k_j}(\mathbf x)) - \min u_\rho ,$$
assures
$$ \max_{x_0 = s, \; \bar x_0 = s'} [u_\rho(\mathbf x) - u_\rho(\bar {\mathbf x})] \leq
k_j\, \text{H\"old}_\theta(A) + 2 \frac{\text{H\"old}_\theta(A)} {1 - \lambda^\theta}\lambda^\theta. $$
Hence, when  $ K = \max k_j $, it follows
$$ \max_{\mathbf x, \bar {\mathbf x} \in \Sigma} [u_\rho(\mathbf x) - u_\rho(\bar {\mathbf x})] \leq
\left(K + \frac{2\lambda^\theta}{1 - \lambda^\theta} \right)\text{H\"old}_\theta(A), $$
that is, the family $ \{u_\rho \} $ has uniformly bounded oscillation.

From the properties demonstrated, we immediately obtain that the family $ \{u_\rho - \max u_\rho \} $
is equicontinuous and uniformly bounded. Note also that
$ u_\rho - \max u_\rho = (\rho - 1) \max u_\rho + \mathcal L_\rho(u_\rho - \max u_\rho) $.
Then, if the function $ u $ (necessarily $\theta$-H\"older) is an accumulation point of $ \{u_\rho - \max u_\rho \} $ when
$ \rho $ tends to 1, we have $ u = a + \mathcal L_1(u) $ for some constant $ a \in \mathbb R $.

It remains to show that $ a = \beta_A $. Put $ \widetilde A = A + u \circ \pi_1 - u \circ \pi_1 \circ \hat \sigma^{-1} $.
Since $ \widetilde A \leq a $, for all $ \hat \mu \in \mathcal M_0 $, we verify
$$ \int_{\hat \Sigma} A(\mathbf y, \mathbf x) \; d\hat\mu(\mathbf y, \mathbf x) =
\int_{\hat \Sigma} \widetilde A(\mathbf y, \mathbf x) \; d\hat\mu(\mathbf y, \mathbf x) \leq a, $$
hence $ \beta_A \leq a $. Besides, observe that
$$ a = \max_{\mathbf y \in \Sigma_{\mathbf x}^*} \widetilde A(\mathbf y, \mathbf x) \; \; \; \; \forall \; \mathbf x \in \Sigma. $$
Thus, given $ \mathbf x^0 \in \Sigma $, take $ \mathbf y^0 \in \Sigma_{\mathbf x^0}^* $ satisfying
$ \widetilde A(\mathbf y^0, \mathbf x^0) = a $. Putting $ \mathbf x^j = \tau_{\mathbf y^{j - 1}}(\mathbf x^{j - 1}) $,
inductively consider $ \mathbf y^j \in \Sigma_{\mathbf x^j}^* $ such that $ \widetilde A(\mathbf y^j, \mathbf x^j) = a $.
Let $ \hat \mu \in \mathcal M $ be an accumulation point of the sequence of probabilities
$$ \hat \mu_k = \frac{1}{k} \sum_{j = 0}^{k - 1} \delta_{(\mathbf y^j, \mathbf x^j)}. $$
Clearly it is true that $ {\displaystyle \int_{\hat \Sigma} \widetilde A(\mathbf y, \mathbf x) \; d\hat\mu(\mathbf y, \mathbf x)} = a $.
Therefore, if we prove that $ \hat \mu \in \mathcal M_0 $, we will obtain $ a \leq \beta_A $. For any $ f \in C^0(\Sigma) $, note then
\begin{eqnarray*}
\left | \int_{\hat \Sigma} \left [ f(\tau_{\mathbf y}(\mathbf x)) - f(\mathbf x) \right ] \; d\hat\mu_k(\mathbf y, \mathbf x) \right |
& = &
\frac{1}{k} \left | \sum_{j = 0}^{k - 1} \left [ f(\tau_{\mathbf y^j}(\mathbf x^j)) - f(\mathbf x^j) \right ] \right |\\
& = &
\frac{1}{k} \left | f(\mathbf x^k) - f(\mathbf x^0) \right | \le \frac{2}{k} \| f \|_0,
\end{eqnarray*}
Now taking the limit when $ k $ tends to infinite, we assure $ \hat \mu \in \mathcal M_0 $ and this finishes the proof.
\end{proof}

The previous result implies the existence of a calibrated sub-action $ u $ for
the forward shift setting \cite{Bousch1, CLT, Jenkinson1}. Indeed, supposing
$ A \in C^\theta(\Sigma) $, observe that we have $ A \circ \tau \in C^\theta(\hat \Sigma) $.
Hence, under the transitivity hypothesis, there exists a function $ u \in C^\theta(\Sigma) $ satisfying
$$ u(\mathbf x) = \min_{\mathbf y \in \Sigma_{\mathbf x}^*} [u(\tau_{\mathbf y}(\mathbf x)) - A \circ \tau (\mathbf y, \mathbf x) + \beta_{A \circ \tau}]. $$
Once $ \beta_{A \circ \tau} = \beta_A = {\displaystyle \max_{\mu \in \mathcal M_\sigma} \int_\Sigma A(\mathbf x) \; d\mu(\mathbf x)} $,
taking $ \mathbf z = \tau_{\mathbf y}(\mathbf x) $, we obtain the
usual expression (see for instance \cite{CLT})
$$ u(\mathbf x) = \min_{\sigma(\mathbf z) = \mathbf x} (u - A + \beta_A)(\mathbf z). $$

The calibrated sub-action notion is an important concept also in relative maximization.
In particular, theorem~\ref{existenciacalibrada} assures a version for the holonomic setting of theorem 17 in \cite{GL}.
Such version will point out that the differential of an alpha application dictates the
asymptotic behavior of the optimal trajectories. We will state the precise result.

We start considering the Fenchel transform of the previous beta function $ \beta_{A, \varphi} $.
Called an alpha application, such function $ \alpha_{A, \varphi}: \mathbb R^n \to \mathbb R $ is defined simply by
$$ \alpha_{A, \varphi}(c) = \min_{h \in \varphi_*(\mathcal M_0)} [\langle c, h \rangle - \beta_{A, \varphi}(h)]. $$
If $ u \in C^0(\Sigma) $ is a calibrated
sub-action, we say that a sequence $ \{\mathbf y^j, \mathbf x^j \} \subset \hat \Sigma $
is an optimal trajectory (associated to the potential $ A $) in the case
$ \mathbf x^j = \tau_{\mathbf y^{j - 1}} (\mathbf x^{j - 1}) $ and
$ u(\mathbf x^j) = u(\mathbf x^{j + 1}) - A(\mathbf y^j, \mathbf x^j) + \beta_A $.
Since the equality $ \alpha_{A, \varphi} (c) = - \beta_{A - \langle c, \varphi \rangle} $ is true,
we can adapt the proof of theorem 17 in \cite{GL} to the present case.
Therefore, under the transitivity  hypothesis, if the potential $ A $ and the constraint $ \varphi $
are H\"older, every optimal trajectory $ \{\mathbf y^j, \mathbf x^j \} $ associated to $ A - \langle c, \varphi \rangle $ satisfies
$$ \lim_{k \to \infty} \frac{1}{k} \sum_{j = 0}^{k - 1} \varphi(\mathbf y^j, \mathbf x^j) = D\alpha_{A, \varphi}(c), $$
in the case the function $ \alpha_{A, \varphi} $ is differentiable at the point $ c \in \mathbb R^n $.

Concluding this section, we would like to say a few words about a
version of Liv\v{s}ic's theorem for the model $ (\hat \Sigma, \hat
\sigma, \mathcal M_0) $.  We will say that a function $ A \in
C^0(\hat \Sigma) $ is cohomologous to a constant $ a \in \mathbb R
$ if there exists a function $ u \in C^0(\Sigma) $ such that
$$ A + u \circ \pi_1 - u \circ \pi_1 \circ \hat \sigma^{-1} = a. $$

\begin{4}
Assume $ \sigma: \Sigma \to \Sigma $ is a transitive subshift of
finite type and suppose that $ A $ is a $\theta$-H\"older
function. Then, $ \text{\Large $\mathit m $}_A = \mathcal M_0 $ if, and only if, $ A $ is cohomologous to
$ \beta_A $.
\end{4}

\begin{proof}
The sufficiency is obvious. Reciprocally, as $ \text{\Large $\mathit m $}_A = \mathcal M_0 $ implies
$ \beta_A = - \beta_{-A} $, consider functions $ u, u' \in C^0(\Sigma) $ satisfying
$$ A + u \circ \pi_1 - u \circ \pi_1 \circ \hat \sigma^{-1} \le \beta_A \; \text { and } \;
\beta_A \le A - u' \circ \pi_1 + u' \circ \pi_1 \circ \hat \sigma^{-1}. $$
Therefore, we have $ (u + u') \circ \pi_1 \le (u + u') \circ \pi_1 \circ \hat \sigma^{-1} $.
In this case, however, the transitivity hypothesis implies that the function $ u + u' $ is identically
equal to a constant $ b $. Since $ u = b - u'  $, from the two above inequalities, it follows that
the potential $ A $ is cohomologous to $ \beta_A $ via the function $ u $.
\end{proof}

\newpage

\textbf{{\large 4. Calibrated Sub-actions and Ma\~n\'e potential}}

\vspace {.5cm}

Using the Ma\~n\'e potential and the set of non-wandering points, we will be able to introduce a family of H\"older calibrated
sub-actions. In the final section, this family will play a crucial role in the classification theorem of calibrated
sub-actions.

\begin{definition}
Given $ \epsilon > 0 $ and  $ \mathbf x,
\bar{\mathbf x} \in \Sigma $, we will call a path
beginning within $\epsilon$ of $ \mathbf x $ and ending at 
$ \bar{\mathbf x} $ an ordered sequence of points
$$ (\mathbf y^0, \mathbf x^0), \ldots, (\mathbf y^{k - 1}, \mathbf x^{k - 1}) \in \hat \Sigma $$ 
satisfying $ \mathbf x^0 = \bar{\mathbf x} $, $ \mathbf x^{j + 1} = \tau_{\mathbf y^j}(\mathbf x^j) $ and 
$ d(\tau_{\mathbf y^{k - 1}} (\mathbf x^{k - 1}), \mathbf x) < \epsilon $. We will denote by 
$ \mathcal P(\mathbf x, \bar{\mathbf x}, \epsilon) $ the set of such paths.
\end{definition}




\begin{definition}
Following \cite{CLT}, a point $ \mathbf x \in
\Sigma $ will be called non-wandering with respect to the
potential $ A \in C^0(\hat \Sigma) $ when, for all $ \epsilon > 0
$, we can determine a path $ \{(\mathbf y^0, \mathbf x^0), \ldots,
(\mathbf y^{k - 1}, \mathbf x^{k - 1}) \} \in \mathcal P(\mathbf
x, \mathbf x, \epsilon) $ such that
$$ \left | \sum_{j = 0}^{k - 1} (A - \beta_A)(\mathbf y^j, \mathbf x^j) \right | < \epsilon. $$
We will denote by $ \Omega(A) $ the set of non-wandering points with respect to $ A $.
\end{definition}

When the potential is H\"older, it is not difficult to see that $ \Omega(A) $ is a compact invariant set.
We will show that such set is indeed not empty.

\begin{3}\label{pontonaoerrante}
If $\sigma: \Sigma \to \Sigma $ is a transitive subshift of finite type,
for any potential $ A \in C^\theta(\hat \Sigma) $, we have $ \Omega(A) \ne \emptyset $.
\end{3}

\begin{proof}
Let $ u \in C^0(\Sigma) $ be a calibrated sub-action obtained from
theorem~\ref{existenciacalibrada}. Fix any point $ \mathbf x^0 \in \Sigma $. Take then $ \mathbf y^0 \in \Sigma_{\mathbf x^0} ^ * $
satisfying the identity $ u(\mathbf x^0) = u(\tau_{\mathbf y^0}(\mathbf x^0)) - A(\mathbf y^0, \mathbf x^0) + \beta_A $.
Denote $ \mathbf x^{j + 1} = \tau_{\mathbf y^j}(\mathbf x^j) $ and proceed in an inductive way determining
a point $ \mathbf y^{j + 1} \in \Sigma_{\mathbf x^{j + 1}} ^ * $ such that
$ u(\mathbf x^{j + 1}) = u(\tau_{\mathbf y^{j + 1}}(\mathbf x^{j + 1})) - A(\mathbf y^{j + 1}, \mathbf x^{j + 1}) + \beta_A $.
Let $ \mathbf x \in \Sigma $ be a limit of some subsequence $ \{\mathbf x^{j_m} \} $.

We claim that $ \mathbf x \in \Omega(A) $. First note that, if $ m_2 > m_1 $,
from the definition of the sequence $ \{\mathbf x^j \} $, we obtain
$$ - \sum_{j = j_{m_1}} ^ {j_{m_2} - 1} (A - \beta_A)(\mathbf y^j, \mathbf x^j) = u(\mathbf x^{j_{m_1}}) - u(\mathbf x^{j_{m_2}}). $$
For a fixed $ \epsilon > 0 $, consider an integer $ l > 0 $ such
that, if $ \mathbf x', \mathbf x '' \in \Sigma $ and $ d(\mathbf x', \mathbf x '') < \lambda^l $, then
$ | u(\mathbf x') - u(\mathbf x'') | < \epsilon/2 $. We can suppose $ l $ is sufficiently large in such way that
$$ \max \left \{\lambda^l, \frac{\text{H\"old}_\theta(A)} {1 - \lambda^\theta} \lambda^{\theta l} \right \} < \frac{\epsilon}{2}. $$
Now take an integer $ m_0 $ sufficiently large such that $ d(\mathbf x^{j_m}, \mathbf x) < \lambda^l/2 $ for all $ m > m_0 $.
Considering integers $ m_2 > m_1 > m_0 $, put $ k = j_{m_2} - j_{m_1} $. Since $ \Sigma_{\mathbf x}^* = \Sigma_{\mathbf x^{j_{m_1}}}^* $,
we choose $ \bar{\mathbf y}^j = \mathbf y^{j_{m_1} + j} $ for $ 0 \le j \le k - 1 $.
Finally, denote $ \bar{\mathbf x}^0 = \mathbf x $ and $ \bar{\mathbf x}^{j + 1} = \tau_{\bar{\mathbf y}^j}(\bar{\mathbf x}^j) $.
Once
$$ d(\tau_{\bar{\mathbf y}^{k - 1}} (\bar{\mathbf x}^{k - 1}), \mathbf x) \le
d(\tau_{\bar{\mathbf y}^{k - 1}} (\bar{\mathbf x}^{k - 1}), \mathbf x^{j_{m_2}}) + d(\mathbf x^{j_{m_2}}, \mathbf x) < \lambda^{k + l} + \lambda^l < \epsilon, $$
it follows that
$ \{(\bar{\mathbf y}^0, \bar{\mathbf x}^0), \ldots, (\bar{\mathbf y}^{k - 1}, \bar{\mathbf x} ^ {k - 1}) \} \in \mathcal P(\mathbf x, \mathbf x, \epsilon) $.
Moreover, since $ d(\mathbf x^{j_{m_1}}, \mathbf x^{j_{m_2}}) < \lambda^l $, we get
\begin{flushleft}
$ {\displaystyle \left | \sum_{j = 0}^{k - 1} (A - \beta_A)(\bar{\mathbf y}^j, \bar{\mathbf x}^j) \right | \leq} $\\
\vspace {.1cm} \centering {$ {\displaystyle \leq \left | \sum_{j = 0}^{k - 1} A(\bar{\mathbf y}^j, \bar{\mathbf x}^j) -
\sum_{j = j_{m_1}} ^ {j_{m_2} - 1} A(\mathbf y^j, \mathbf x^j) \right | + | u(\mathbf x^{j_{m_1}}) - u(\mathbf x^{j_{m_2}}) | <} $} \\
\vspace {.1cm} \raggedleft {$ {\displaystyle < \frac{\text{H\"old}_\theta(A)} {1 - \lambda^\theta} \lambda^{\theta l} + \frac{\epsilon}{2} < \epsilon} $.}
\end{flushleft}
Therefore, $ \mathbf x \in \Omega(A) $.
\end{proof}

The following definition is also inspired in \cite{CLT}.

\begin{definition}
We call Ma\~n\'e potential the function $ S_A: \Sigma \times \Sigma \to \mathbb R \cup \{ \pm \infty \} $
defined by
$$ S_A(\mathbf x, \bar{\mathbf x}) = \lim_{\epsilon \to 0} S_A^\epsilon(\mathbf x, \bar{\mathbf x}), $$
where
$$ S_A^\epsilon(\mathbf x, \bar{\mathbf x}) =
\inf_{\{(\mathbf y^0, \mathbf x^0), \ldots, (\mathbf y^{k - 1}, \mathbf x^{k - 1}) \} \in \mathcal P(\mathbf x, \bar{\mathbf x}, \epsilon)}
\left [- \sum_{j = 0}^{k - 1} (A - \beta_A)(\mathbf y^j, \mathbf x^j) \right]. $$
\end{definition}

Note that $ \Omega(A) = \{\mathbf x \in \Sigma: S_A(\mathbf x, \mathbf x) = 0 \} $.

As we will see soon the Ma\~n\'e potential will provide, for a H\"older potential, a one-parameter family of
equally H\"older sub-actions. Before that we need some properties.

Let $ u \in C^0(\Sigma) $ be a sub-action for the potential $ A \in C^0(\hat \Sigma) $. We say
that the point $ \mathbf x \in \Sigma $ is $u$-connected to the
point $ \bar{\mathbf x} \in \Sigma $, and we indicate this by
$ \mathbf x \stackrel{u}{\rightarrow} \bar{\mathbf x} $, when, for
every $ \epsilon > 0 $, we can determine a path
$ \{(\mathbf y^0, \mathbf x^0), \ldots, (\mathbf y^{k - 1}, \mathbf x^{k - 1}) \} \in \mathcal P(\mathbf x, \bar{\mathbf x}, \epsilon) $ such that
$$ \left | \sum_{j = 0}^{k - 1} (A - \beta_A)(\mathbf y^j, \mathbf x^j) - (u(\mathbf x) - u(\bar{\mathbf x})) \right | < \epsilon. $$
Note that $ \mathbf x \in \Omega(A) $ implies $ \mathbf x \stackrel{u}{\rightarrow} \mathbf x $ for any sub-action $ u $.

\begin{3}\label{desigualdadepotencialsubacao}
Let $ u \in C^0(\Sigma) $ be a sub-action for a potential $ A \in C^0(\hat \Sigma) $.
Then, for any $ \mathbf x, \bar{\mathbf x} \in \Sigma $, we have
$ S_A(\mathbf x, \bar{\mathbf x}) \ge u(\bar{\mathbf x}) - u(\mathbf x) $.
Moreover, the equality is true if, and only if, $ \mathbf x \stackrel{u}{\rightarrow} \bar{\mathbf x} $.
\end{3}

Before the proof of this lemma, we would like just to point out another important property of Ma\~n\'e potential: if $ A $ is a $\theta$-H\"older potential, then
$ S_A(\mathbf x, \bar{\bar{\mathbf x}}) \leq S_A(\mathbf x, \bar{\mathbf x}) + S_A(\bar{\mathbf x}, \bar{\bar{\mathbf x}}) $
for any points $ \mathbf x, \bar{\mathbf x}, \bar{\bar{\mathbf x}} \in \Sigma $. We leave for the reader the demonstration of this simple fact.

\begin{proof}
Fix $ \rho > 0 $. Take $ \epsilon \in (0, \rho) $ such that $ | u(\mathbf x') - u(\mathbf x '') | < \rho $, when
$ \mathbf x', \mathbf x '' \in \Sigma $ satisfy $ d(\mathbf x', \mathbf x '') < \epsilon $.
Consider now any path
$$ \{(\mathbf y^0, \mathbf x^0), \ldots, (\mathbf y^{k - 1}, \mathbf x^{k - 1}) \} \in \mathcal P(\mathbf x, \bar{\mathbf x}, \epsilon). $$
Once
$$ u(\bar{\mathbf x}) - u(\mathbf x) - \rho < u(\mathbf x^0) - u(\tau_{\mathbf y^{k - 1}} (\mathbf x^{k - 1})) \le
- \sum_{j = 0}^{k - 1} (A - \beta_A)(\mathbf y^j, \mathbf x^j), $$
it follows that $ u(\bar{\mathbf x}) - u(\mathbf x) - \rho \le S_A(\mathbf x, \bar{\mathbf x}) $. Taking
$ \rho $ arbitrarily small, we obtain the inequality of the lemma.

If $ S_A(\mathbf x, \bar{\mathbf x}) = u(\bar{\mathbf x}) - u(\mathbf x) $, from the definition of the Ma\~n\'e potential,
immediately we get $ \mathbf x \stackrel{u}{\rightarrow} \bar{\mathbf x} $. Reciprocally, suppose that $ \mathbf x $ is
$u$-connected to $ \bar{\mathbf x} $. Take then $ \rho > 0 $. Given $ \epsilon \in (0, \rho) $, we can choose a path
$$ \{(\mathbf y^0, \mathbf x^0), \ldots, (\mathbf y^{k - 1}, \mathbf x^{k - 1}) \} \in \mathcal P(\mathbf x, \bar{\mathbf x}, \epsilon) $$
satisfying
$$ \left | \sum_{j = 0}^{k - 1} (A - \beta_A)(\mathbf y^j, \mathbf x^j) - (u(\mathbf x) - u(\bar{\mathbf x})) \right | < \epsilon. $$
Observe that
$$ - \sum_{j = 0}^{k - 1} (A - \beta_A)(\mathbf y^j, \mathbf x^j) < u(\bar{\mathbf x}) - u(\mathbf x) + \epsilon < u(\bar{\mathbf x}) - u(\mathbf x) + \rho. $$
Thus, we verify $ S_A(\mathbf x, \bar{\mathbf x}) \leq u(\bar{\mathbf x}) - u(\mathbf x) + \rho $.
As $ \rho $ can be taken arbitrarily small, we finally get the equality claimed by the lemma.
\end{proof}

We present now the main result of this section.

\begin{4}\label{familiacalibradas}
Suppose $ \sigma: \Sigma \to \Sigma $ is a
transitive subshift of finite type. Let $ A $ be a
$\theta$-H\"older potential. Then, for each $\mathbf x \in
\Omega(A) $, the function $S_A(\mathbf x, \cdot) $  is a
$\theta$-H\"older calibrated sub-action.
\end{4}

\begin{proof}
Fix a point $ \mathbf x \in \Omega(A) $. We must show first
that $ S_A(\mathbf x, \cdot) $ is a well defined real function. Thanks to lemma~\ref{desigualdadepotencialsubacao}, 
we only need to assure that $ S_A(\mathbf x, \bar{\mathbf x}) < + \infty $ for any $ \bar{\mathbf x} \in \Sigma $.

Take $ \epsilon > 0 $ arbitrary. For a fixed value $ \epsilon' \in (0, \lambda] $, consider a path
$ \{(\mathbf y^0, \mathbf x^0), \ldots, (\mathbf y^{k - 1}, \mathbf x^{k - 1}) \} \in \mathcal P(\mathbf x, \bar{\mathbf x}, \epsilon') $
satisfying
$$ - \sum_{j = 0}^{k - 1} (A - \beta_A)(\mathbf y^j, \mathbf x^j) < S_A^{\epsilon'} (\mathbf x, \bar{\mathbf x}) + \epsilon. $$
As $ \mathbf x \in \Omega(A) $, we can take
$ \{(\bar{\mathbf y}^0, \bar{\mathbf x}^0), \ldots, (\bar{\mathbf y}^{\bar k - 1}, \bar{\mathbf x}^{\bar k - 1}) \} \in
\mathcal P(\mathbf x, \mathbf x, \epsilon/2) $, with $ \lambda^{\bar k} \epsilon' < \epsilon/2 $, such that
$$ \left | \sum_{j = 0}^{\bar k - 1} (A - \beta_A)(\bar{\mathbf y}^j, \bar{\mathbf x}^j) \right | < \frac{\epsilon}{2}. $$

Thus, we define $ \mathbf y^j = \bar{\mathbf y}^{j - k} $ for $ k \le j < k + \bar k $.
Observe that we have $ \mathbf y^k = \bar{\mathbf y}^0 \in \Sigma_{\bar{\mathbf x}^0}^* = \Sigma_{\tau_{\mathbf y^{k - 1}} (\mathbf x^{k - 1})}^* $.
Therefore, we can put $ \mathbf x^{j + 1} = \tau_{\mathbf y^j}(\mathbf x^j) $ for $ k - 1 \le j < k + \bar k - 1 $.

We claim that $\{(\mathbf y^0, \mathbf x^0), \ldots, (\mathbf y^{k + \bar k - 1}, \mathbf x^{k + \bar k - 1}) \} \in
\mathcal P(\mathbf x, \bar{\mathbf x}, \epsilon) $. Indeed,
\begin{flushleft}
${\displaystyle d(\tau_{\mathbf y^{k + \bar k - 1}} (\mathbf x^{k + \bar k - 1}), \mathbf x) \le} $ \\
\vspace {.1cm}
\centering {${\displaystyle \le
d(\tau_{\mathbf y^{k + \bar k - 1}} (\mathbf x^{k + \bar k - 1}), \tau_{\bar{\mathbf y}^{\bar k - 1}} (\bar{\mathbf x}^{\bar k - 1})) +
d(\tau_{\bar{\mathbf y}^{\bar k - 1}} (\bar{\mathbf x}^{\bar k - 1}), \mathbf x) <} $} \\
\vspace {.1cm}
\raggedleft {${\displaystyle < \lambda^{\bar k} \epsilon' + \frac{\epsilon}{2} < \epsilon} $.}
\end{flushleft}
Besides, without difficulty we verify
$$ \left | \sum_{j = k}^{k
+ \bar k - 1} A(\mathbf y^j, \mathbf x^j) - \sum_{j = 0}^{\bar k -
1} A(\bar{\mathbf y}^j, \bar{\mathbf x}^j) \right | \leq
\frac{\text{H\"old}_\theta(A)} {1 - \lambda^\theta}(\epsilon') ^
\theta. $$

Hence, we immediately have
$$ S_A^\epsilon(\mathbf x, \bar{\mathbf x}) \le - \sum_{j = 0}^{k + \bar k - 1} (A - \beta_A)(\mathbf y^j, \mathbf x^j) <
\frac{\text{H\"old}_\theta(A)} {1 - \lambda^\theta}(\epsilon') ^ \theta + S_A^{\epsilon'} (\mathbf x, \bar{\mathbf x}) + \frac{3}{2}\epsilon, $$
which yields
$$ S_A(\mathbf x, \bar{\mathbf x}) \le \frac{\text{H\"old}_\theta(A)} {1 - \lambda^\theta}(\epsilon')^\theta + S_A^{\epsilon'} (\mathbf x, \bar{\mathbf x}). $$
As the right hand side is finite, the application $ S_A(\mathbf x, \cdot) $ is well defined.

We claim that it is indeed a $\theta$-H\"older function.
Take points $ \bar{\mathbf x}, \bar{\bar{\mathbf x}} \in \Sigma $ such
that $ d(\bar{\mathbf x}, \bar{\bar{\mathbf x}}) \le \lambda $.
Consider a fixed $ \rho > 0 $. Given $ \epsilon > 0 $, we can find a
path $ \{(\mathbf y^0, \mathbf x^0), \ldots, (\mathbf y^{k - 1}, \mathbf x^{k - 1}) \} \in \mathcal P(\mathbf x, \bar{\mathbf x}, \epsilon) $,
with $ \lambda^{k + 1} < \epsilon $, such that
$$ - \sum_{j = 0}^{k - 1} (A - \beta_A)(\mathbf y^j, \mathbf x^j) < S_A^{\epsilon}(\mathbf x, \bar{\mathbf x}) + \rho. $$

Taking $ \bar{\mathbf y}^j = \mathbf y^j $ for $ 0 \le j < k $, we
write $ \bar{\mathbf x}^0 = \bar{\bar{\mathbf x}} $ and, finally,
we define $ \bar{\mathbf x}^{j + 1} = \tau_{\bar{\mathbf
y}^j}(\bar{\mathbf x}^j) $ when $ 0 \le j < k - 1 $. It is easy to
confirm that $ \{(\bar{\mathbf y}^0,\bar{\mathbf x}^0), \ldots,
(\bar{\mathbf y}^{k - 1},\bar{\mathbf x}^{k - 1}) \} \in \mathcal
P(\mathbf x, \bar{\bar{\mathbf x}}, 2\epsilon) $, as well as
$$ - \sum_{j = 0}^{k - 1} A(\mathbf y^j, \mathbf x^j) \ge - \sum_{j = 0}^{k - 1} A(\bar{\mathbf y}^j, \bar{\mathbf x}^j) -
\frac{\text{H\"old}_\theta(A)}{1 - \lambda^\theta} d(\bar{\mathbf x}, \bar{\bar{\mathbf x}})^\theta. $$

Therefore, we verify the following inequalities
\begin{eqnarray*}
S_A(\mathbf x, \bar{\mathbf x}) & \ge & S_A^\epsilon(\mathbf x, \bar{\mathbf x}) \\
& > & - \sum_{j = 0}^{k - 1} (A - \beta_A)(\mathbf y^j, \mathbf x^j) - \rho \\
& \ge & - \sum_{j = 0}^{k - 1} (A - \beta_A)(\bar{\mathbf y}^j, \bar{\mathbf x}^j) -
\frac{\text{H\"old}_\theta(A)}{1 - \lambda^\theta} d(\bar{\mathbf x}, \bar{\bar{\mathbf x}})^\theta - \rho \\
& \ge & S_A^{2\epsilon}(\mathbf x, \bar{\bar{\mathbf x}}) -
\frac{\text{H\"old}_\theta(A)}{1 - \lambda^\theta} d(\bar{\mathbf x}, \bar{\bar{\mathbf x}})^\theta - \rho.
\end{eqnarray*}

Since $ \epsilon $ and $ \rho $ can be considered (in such order) arbitrarily small, we get
$$ S_A(\mathbf x, \bar{\mathbf x}) - S_A(\mathbf x, \bar{\bar{\mathbf x}}) \ge
- \frac{\text{H\"old}_\theta(A)}{1 - \lambda^\theta} d(\bar{\mathbf x}, \bar{\bar{\mathbf x}})^\theta. $$
It follows at once that $ S_A(\mathbf x, \cdot) \in C^\theta(\Sigma) $.

It remains to show that the application $ S_A(\mathbf x, \cdot) $ is a calibrated sub-action.

Fix a point $ (\bar{\mathbf y}, \bar{\mathbf x}) \in \hat \Sigma $. When
$ \{(\mathbf y^1, \mathbf x^1), \ldots, (\mathbf y^k, \mathbf x^k) \} \in \mathcal P(\mathbf x, \tau_{\bar{\mathbf y}} (\bar{\mathbf x}), \epsilon) $,
put $ \mathbf y^0 = \bar{\mathbf y} $, $ \mathbf x^0 = \bar{\mathbf x} $. We point out that
\begin{eqnarray*}
A(\bar{\mathbf y}, \bar{\mathbf x}) - \beta_A & = &
\sum_{j = 0}^k (A - \beta_A)(\mathbf y^j, \mathbf x^j) - \sum_{j = 0}^{k - 1} (A - \beta_A)(\mathbf y^{j + 1}, \mathbf x^{j + 1}) \\
& \le &  - \sum_{j = 0}^{k - 1} (A - \beta_A)(\mathbf y^{j + 1}, \mathbf x^{j + 1}) - S_A^\epsilon(\mathbf x, \bar{\mathbf x}).
\end{eqnarray*}
As the path is arbitrary, we have $ A(\bar{\mathbf y}, \bar{\mathbf x}) - \beta_A \le
S_A^\epsilon(\mathbf x, \tau_{\bar{\mathbf y}} (\bar{\mathbf x})) - S_A^\epsilon(\mathbf x, \bar{\mathbf x}) $.
Hence, taking limit, we show that $ S_A(\mathbf x, \cdot) $ is indeed a sub-action for the potential $ A $.

In order to verify that it is a calibrated sub-action, we should be able to
determine, for each $ \bar{\mathbf x} \in \Sigma $, a point $ \bar{\mathbf y} \in \Sigma_{\bar{\mathbf x}}^* $ accomplishing the equality
$ S_A(\mathbf x, \bar{\mathbf x}) = S_A(\mathbf x, \tau_{\bar{\mathbf y}} (\bar{\mathbf x})) - A(\bar{\mathbf y}, \bar{\mathbf x}) + \beta_A $.
Given $ \epsilon > 0 $, consider a path
$ \{(\mathbf y^0, \mathbf x^0), \ldots, (\mathbf y^{k - 1}, \mathbf x^{k - 1}) \} \in \mathcal P(\mathbf x, \bar{\mathbf x}, \epsilon) $
 such that
$$ - \sum_{j = 0}^{k - 1} (A - \beta_A)(\mathbf y^j, \mathbf x^j) < S_A^\epsilon(\mathbf x, \bar{\mathbf x}) + \epsilon. $$
This defines a family $\{\mathbf y^0 \} _{\epsilon > 0} \subset \Sigma_{\bar{\mathbf x}} ^ * $. Take $ \bar{\mathbf y} \in \Sigma_{\bar{\mathbf x}}^* $
an accumulation point of this family when $ \epsilon $ tends to 0. Observe that
$$ S_A^\epsilon(\mathbf x, \tau_{\mathbf y^0}(\bar{\mathbf x})) - (A - \beta_A)(\mathbf y^0, \bar{\mathbf x}) \le
- \sum_{j = 0}^{k - 1} (A - \beta_A)(\mathbf y^j, \mathbf x^j). $$
As $ \tau_{\mathbf y^0}(\bar{\mathbf x}) = \tau_{\bar{\mathbf y}} (\bar{\mathbf x}) $ for $ \epsilon $ sufficiently small, we can focus on
$$ S_A^\epsilon(\mathbf x, \tau_{\bar{\mathbf y}} (\bar{\mathbf x})) - (A - \beta_A)(\mathbf y^0, \bar{\mathbf x}) <
S_A^\epsilon(\mathbf x, \bar{\mathbf x}) + \epsilon . $$
So taking $ \epsilon $ arbitrarily small, we finish the proof.
\end{proof}

\vspace {1cm}

\textbf{{\large 5. Sub-actions and Supports}}

\vspace {.5cm}

This section is dedicated to the analysis of relationships between
sub-actions and supports of holonomic probabilities.
An unifying element of these concepts continues to be the contact
locus notion.

\begin{definition}
Given a sub-action $ u \in C^0(\Sigma) $ for a
potential $ A \in C^0(\hat \Sigma) $, consider the function $ A^u
= A + u \circ \pi_1 - u \circ \pi_1 \circ \hat \sigma^{-1} $. We
call the set $ \mathbb M_A(u) = (A^u)^{-1}(\beta_A) $ the contact
locus of the sub-action $ u $.
\end{definition}

The contact locus is just the set where the usual inequality defining a sub-action
becomes an equality. It plays an important role in the localization of the
support of maximing holonomic probabilities.

\begin{4}\label{contatomaximizantes}
If $ u \in C^0(\Sigma) $ is a sub-action for a potential $ A \in C^0(\hat \Sigma) $, then
$$ \text{\Large $\mathit m $}_A = \left \{\hat \mu \in \mathcal M_0: \text{supp}(\hat \mu) \subset \mathbb M_A(u) \right \}. $$
\end{4}

The proof of this statement is reduced to the well known fact
according to which is zero almost everywhere a measurable non negative
function whose integral is zero.

We aim now a classification theorem for calibrated sub-actions.
We start presenting a result which supplies a representation
formula for these sub-actions.

\begin{2}\label{formularepresentacao}
If $ u \in C^0(\Sigma) $ is a calibrated sub-action for a $\theta$-H\"older potential $ A $, then
$$ u(\bar{\mathbf x}) = \inf_{\mathbf x \in \Omega(A)} [u(\mathbf x) + S_A(\mathbf x, \bar{\mathbf x})]. $$
\end{2}

\begin{proof}
Thanks to lemma~\ref{desigualdadepotencialsubacao}, it immediately follows that
$$ u(\bar{\mathbf x}) \le \inf_{\mathbf x \in \Omega(A)} [u(\mathbf x) + S_A(\mathbf x, \bar{\mathbf x})]. $$
Besides, the identity will be true if there exists a point $ \mathbf x \in \Omega(A) $ satisfying
$ \mathbf x \stackrel{u}{\rightarrow} \bar{\mathbf x} $.

Consider $ \{(\mathbf y^j, \mathbf x^j) \} \subset \hat \Sigma $ an
optimal trajectory associated to the potential $ A $ such that $ \mathbf x^0 = \bar{\mathbf x} $.
Denote by $ \mathbf x \in \Sigma $ the limit of a subsequence $ \{\mathbf x^{j_m} \} $.

Lemma~\ref{pontonaoerrante} shows that $ \mathbf x \in \Omega(A) $. So we only have to
prove that $ \mathbf x \stackrel{u}{\rightarrow} \bar{\mathbf
x} $. Fix $ \epsilon > 0 $ and  choose an integer $ l > 0 $ in such
way that $ | u(\mathbf x') - u(\mathbf x '') | < \epsilon $ when
$\mathbf x', \mathbf x '' \in \Sigma $ satisfy $ d(\mathbf x', \mathbf x '') < \lambda^l $.
Assume $ l $ also accomplishes $ \lambda^l < \epsilon $.
Take $ m $ sufficiently large such that $ d(\mathbf x^{j_m}, \mathbf x) < \lambda^l $. Put
$ k = j_m $.

Observe that $ d(\tau_{\mathbf y^{k - 1}} (\mathbf x^{k - 1}), \mathbf x) = d(\mathbf x^{j_m}, \mathbf x) < \epsilon $.
Therefore, we assure
$ \{(\mathbf y^0, \mathbf x^0), \ldots, (\mathbf y^{k - 1}, \mathbf x^{k - 1}) \} \in \mathcal P(\mathbf x, \bar{\mathbf x}, \epsilon) $.
As
$$\sum_{j = 0}^{k - 1} (A - \beta_A)(\mathbf y^j, \mathbf
x^j) - (u(\mathbf x^k) - u(\bar{\mathbf x})) = 0,
$$
we obtain
$$ \left | \sum_{j = 0}^{k - 1} (A - \beta_A)(\mathbf y^j, \mathbf x^j) - (u(\mathbf x) - u(\bar{\mathbf x})) \right | =
| u(\mathbf x^{j_m}) - u(\mathbf x) | < \epsilon, $$
which finishes the proof.
\end{proof}

The following immediate corollary indicates the importance of the set $ \Omega(A) $ in
the analysis of calibrated sub-actions.

\begin{5}
Let $ u, u' \in C^0(\Sigma) $ be calibrated sub-actions for a potential $ A \in C^\theta(\hat \Sigma) $.
If $ u \le u' $ on $ \Omega(A) $, then $ u \le u' $ everywhere on $ \Sigma $. In particular,
if we have $ u|_{\Omega(A)} = u'|_{\Omega(A)}$, then both sub-actions are equal.
\end{5}

The theorem~\ref{formularepresentacao} admits a reciprocal.

\begin{2}\label{reciprocarepresentacao}
Let $ \sigma: \Sigma \to \Sigma $ be a transitive subshift of finite type. Consider
a potencial $ A \in C^\theta(\hat \Sigma) $. Assume that the function
$ f: \Omega(A) \to \mathbb R $ has a finite lower bound. Then
$$ u(\bar{\mathbf x}) = \inf_{\mathbf x \in \Omega(A)} [f(\mathbf x) + S_A(\mathbf x, \bar{\mathbf x})] $$
defines a $\theta$-H\"older calibrated sub-action. Moreover, if
$ f(\bar{\mathbf x}) - f(\mathbf x) \le S_A(\mathbf x, \bar{\mathbf x}) $ for any
$ \mathbf x, \bar{\mathbf x} \in \Omega(A) $, then $ u = f $ on $ \Omega(A) $.
\end{2}

\begin{proof} The good definition of $ u: \Sigma \to \mathbb R $ is clear. We will show it
is a H\"older function. Fix $ \epsilon > 0 $. Given $ \bar{\mathbf x}, \bar{\bar{\mathbf x}} \in \Sigma $ with
$ d(\bar{\mathbf x}, \bar{\bar{\mathbf x}}) \le \lambda $, take a point
$ \mathbf x \in \Omega(A) $ such that
$ f(\mathbf x) + S_A(\mathbf x, \bar{\bar{\mathbf x}}) < u(\bar{\bar{\mathbf x}}) + \epsilon $.
It follows from the proof of proposition~\ref{familiacalibradas} that
$$ u(\bar{\mathbf x}) - u(\bar{\bar{\mathbf x}}) - \epsilon <
S_A(\mathbf x, \bar{\mathbf x}) - S_A(\mathbf x, \bar{\bar{\mathbf x}}) \le
\frac{\text{H\"old}_\theta(A)} {1 - \lambda^\theta} d(\bar{\mathbf x}, \bar{\bar{\mathbf x}})^\theta. $$
As $ \epsilon $ is arbitrary, we get $ u \in C^\theta(\Sigma) $.

In fact, $ u $ is a sub-action for the potential $ A $.
Consider a point $ (\bar{\mathbf y}, \bar{\mathbf x}) \in \hat \Sigma $
and $ \epsilon > 0 $. Choose $ \mathbf x \in \Omega(A) $ satisfying
$ f(\mathbf x) + S_A(\mathbf x, \tau_{\bar{\mathbf y}} (\bar{\mathbf x})) <
u(\tau_{\bar{\mathbf y}} (\bar{\mathbf x})) + \epsilon $.
Since
$$ u(\bar{\mathbf x}) - u(\tau_{\bar{\mathbf y}} (\bar{\mathbf x})) - \epsilon <
S_A(\mathbf x, \bar{\mathbf x}) - S_A(\mathbf x, \tau_{\bar{\mathbf y}} (\bar{\mathbf x})) \le
\beta_A - A(\bar{\mathbf y}, \bar{\mathbf x}), $$
the claim follows when $ \epsilon $ tends to 0.

The calibrated character of $ u $ is also a consequence of proposition~\ref{familiacalibradas}. Indeed, take
$ \bar{\mathbf x} \in \Sigma $, and choose a point $ \mathbf x^j \in \Omega(A) $ such that
$$ f(\mathbf x^j) + S_A(\mathbf x^j, \bar{\mathbf x}) < u(\bar{\mathbf x}) + \frac{1}{j}. $$
Now, for each index $ j $, take a point $ \mathbf y^j \in \Sigma_{\bar{\mathbf x}} ^ * $ satisfying
$$ S_A(\mathbf x^j, \bar{\mathbf x}) = S_A(\mathbf x^j, \tau_{\mathbf y^j}(\bar{\mathbf x})) - A(\mathbf y^j, \bar{\mathbf x}) + \beta_A. $$
Finally, let $ \bar{\mathbf y} \in \Sigma_{\bar{\mathbf x}}^* $ be an accumulation point of the sequence $ \{\mathbf y^j \} $.
As $ u(\tau_{\mathbf y^j}(\bar{\mathbf x})) \le f(\mathbf x^j) + S_A(\mathbf x^j, \tau_{\mathbf y^j}(\bar{\mathbf x})) $, we verify
$$ u(\tau_{\mathbf y^j}(\bar{\mathbf x})) - A(\mathbf y^j, \bar{\mathbf x}) + \beta_A < u(\bar{\mathbf x}) + \frac{1}{j}. $$
Therefore, $ u(\tau_{\bar{\mathbf y}} (\bar{\mathbf x})) - A(\bar{\mathbf y}, \bar{\mathbf x}) + \beta_A \leq u(\bar{\mathbf x}) $.

At last, suppose that $ f(\bar{\mathbf x}) - f(\mathbf x) \le S_A(\mathbf x, \bar{\mathbf x}) $ for any
$ \mathbf x, \bar{\mathbf x} \in \Omega(A) $. Hence, the inequalities $ u(\bar{\mathbf x}) \le
f(\bar{\mathbf x}) \le f(\mathbf x) + S_A(\mathbf x, \bar{\mathbf x}) $ are valid for all $ \mathbf x \in \Omega(A) $,
which implies immediately $ u = f $ on $ \Omega(A) $.
\end{proof}

One of the main consequences of the previous theorem is a kind of
\emph{H\"older supremacy} for sub-actions that we will state
bellow. This result corresponds to the well known fact in
Lagrangian Aubry-Mather theory according to which a weak KAM solution is
differentiable in the Aubry set (see \cite{CI}).

\begin{5}
Suppose $ \sigma: \Sigma \to \Sigma $ is a
transitive subshift of finite type. If $ u \in C^0(\Sigma) $
is a sub-action for a potential $ A \in C^\theta(\hat \Sigma) $,
then $ u |_{\Omega(A)} $ is $\theta$-H\"older.
\end{5}

Allow us to indicate another immediate consequence of theorem~\ref{reciprocarepresentacao}.

\begin{5}\label{comportamentocalibradageral}
Let $ \sigma: \Sigma \to \Sigma $ be a transitive subshift of finite type. 
Assume $ u \in C^0(\Sigma) $ is a sub-action for a $\theta$-H\"older potential $ A $.
Then, for every point $ \mathbf x \in \Omega(A) $, we verify
$$ u(\mathbf x) = \min_{\mathbf y \in \Sigma_{\mathbf x}^*} [u(\tau_{\mathbf y}(\mathbf x)) - A(\mathbf y, \mathbf x) + \beta_A]. $$
\end{5}

Theorems~\ref{formularepresentacao} and~\ref{reciprocarepresentacao} assure that every calibrated sub-action for a
H\"older potential $A$ is also H\"older. Moreover, we have a complete
description of the set of these sub-actions.

\begin{2}
Consider $ \sigma: \Sigma \to \Sigma $ a
transitive subshift of finite type and $ A : \Sigma \to \mathbb{R}$ a $\theta$-H\"older potential.
Then, there exists a bijective and isometric correspondence between the set of
calibrated sub-actions for $ A $ and the set of functions
$ f \in C^0(\Omega(A)) $ satisfying $ f(\bar{\mathbf x}) - f(\mathbf x) \le S_A(\mathbf x, \bar{\mathbf x}) $,
for all points $ \mathbf x, \bar{\mathbf x} \in \Omega(A) $.
\end{2}

\begin{proof}
Let us analyze the correspondence
$$ f \mapsto u_f = \inf_{\mathbf x \in \Omega(A)} [f(\mathbf x) + S_A(\mathbf x, \cdot)]. $$
It follows from theorem~\ref{reciprocarepresentacao} that such correspondence is well
defined and injective. From theorem~\ref{formularepresentacao} we get that it is
surjective. Besides, the correspondence is an  isometry.
Indeed, fixing $ \epsilon > 0 $, if $ \bar{\mathbf x} \in \Sigma $, take a point
$ \mathbf x \in \Omega(A) $ such that $f(\mathbf x) + S_A(\mathbf x, \bar{\mathbf x}) < u_f(\bar{\mathbf x}) + \epsilon $.
Therefore,
$$ u_g(\bar{\mathbf x}) - u_f(\bar{\mathbf x}) - \epsilon <
g(\mathbf x) - f(\mathbf x) \leq \| f - g \|_0. $$
When $ \epsilon $ tends to 0, since $ \bar{\mathbf x} $ is arbitrary
and since we can interchange the roles of $ f $ and $ g $, we see that
$ \|u_f - u_g\|_0 \le \| f - g\|_0 $. On the other hand, as $ u_f | _{\Omega(A)} = f $
and $ u_g | _{\Omega(A)} = g $, we verify $ \| u_f - u_g\|_0 \ge \| f - g\|_0 $.
\end{proof}

In \cite{Contreras}, Contreras characterizes the weak KAM
solutions of the Hamilton-Jacobi equation in terms of their values
at each static class and the values of the action potential of
Ma\~n\'e. The result we presented above describe similar
property for our holonomic setting.

As announced just before the statement of theorem~\ref{existenciacalibrada}, under the transitive hypothesis, 
there always exists a calibrated sub-action of maximal character for a H\"older potential.
We only need to consider the following one
$$ u_0 = \inf_{\mathbf x \in \Omega(A)} S_A(\mathbf x, \cdot). $$
Indeed, it is clear that $ u_0 \le 0 $ on $ \Omega(A) $. Moreover, if we take any sub-action 
$ u \in C^0(\Sigma) $ satisfying $ u|_{\Omega(A)} \le 0 $, since $ u(\bar{\mathbf x}) \le  
u(\mathbf x) + S_A(\mathbf x, \bar{\mathbf x}) \le S_A(\mathbf x, \bar{\mathbf x}) $ for 
$ \mathbf x \in \Omega(A) $ and $ \bar{\mathbf x} \in \Sigma $, we verify $ u \le u_0 $.

Now we will focus also on the support of maximizing holonomic probabilities in order
to complete our investigation. We need just two lemmas.

\begin{3}\label{suporteholonomicas}
Suppose $ \hat \mu \in \mathcal M_0 $.
Then, almost every point
$ (\mathbf y, \mathbf x) \in \text{supp}(\hat \mu) $ is of the form
$ (\mathbf y, \tau_{\bar{\mathbf y}} (\bar{\mathbf x})) $, with
$ (\bar{\mathbf y}, \bar{\mathbf x}) \in \text{supp}(\hat \mu) $.
\end{3}

\begin{proof}
Consider the set
$$ \hat R = \left \{(\mathbf y, \mathbf x) \in \text{supp}(\hat \mu):
\mathbf x \ne \tau_{\bar{\mathbf y}} (\bar{\mathbf x}) \; \; \; \forall \; (\bar{\mathbf y}, \bar{\mathbf x}) \in \text{supp}(\hat \mu) \right \}. $$
Suppose $ \hat \mu (\hat R) = \epsilon > 0 $. Put $ R = \pi_1 (\hat R) $.
Consider $ D \subset \Sigma $ a compact subset and $ E \subset \Sigma $ an open subset satisfying $ D \subset R \subset E $
with $(\hat \mu \circ \pi_1^{-1}) (E - D) < \epsilon/2 $. Take then a function $ f \in C^0(\Sigma, [0, 1]) $ such that
$ f|_D \equiv 1 $ and $ f|_{\Sigma - E} \equiv 0 $. Once $ \pi_1^{-1}(R) \cap \text{supp}(\hat \mu) = \hat R $, we get
$$ \int_{\hat \Sigma} f(\mathbf x) \; d\hat\mu(\mathbf y, \mathbf x) \ge \hat\mu(\pi_1^{-1}(D)) \ge
\hat\mu(\pi_1^{-1}(R)) - \hat\mu(\pi_1^{-1}(E - D)) > \frac{\epsilon}{2}. $$

Thus, consider a sequence of functions $ \{f_j \} \subset C^0(\Sigma, [0, 1]) $ such that $ f_j \uparrow \chi_{E - D} $.
By the monotonous convergence theorem, we obtain
\begin{eqnarray*}
\int_{\hat \Sigma} \chi_{E - D}(\tau_{\mathbf y}(\mathbf x)) \; d\hat\mu(\mathbf y, \mathbf x)
& = & \lim_{j \to \infty} \int_{\hat \Sigma} f_j(\tau_{\mathbf y}(\mathbf x)) \; d\hat\mu(\mathbf y, \mathbf x) \\
& = & \lim_{j \to \infty} \int_{\hat \Sigma} f_j(\mathbf x) \; d\hat\mu(\mathbf y, \mathbf x) \\
& = & \hat\mu(\pi_1^{-1}(E - D)) < \frac{\epsilon}{2}.
\end{eqnarray*}
Note that, from the definition of $ R $, we have
$ {\displaystyle \int_{\text{supp}(\hat \mu)} \chi_R(\tau_{\mathbf y}(\mathbf x)) \; d\hat\mu(\mathbf y, \mathbf x)} = 0 $.
Hence, as $ 0 \le f \le \chi_E $, we verify
\begin{eqnarray*}
\int_{\hat \Sigma} f(\tau_{\mathbf y}(\mathbf x)) \; d\hat\mu(\mathbf y, \mathbf x) & \le &
\int_{\text{supp}(\hat \mu)} \chi_{E - R}(\tau_{\mathbf y}(\mathbf x)) \; d\hat\mu(\mathbf y, \mathbf x) \\
& \le & \int_{\text{supp}(\hat \mu)} \chi_{E - D}(\tau_{\mathbf y}(\mathbf x)) \; d\hat\mu(\mathbf y, \mathbf x) < \frac{\epsilon}{2}.
\end{eqnarray*}
However, since $ f \in C^0(\Sigma) $ and $ \hat \mu \in \mathcal M_0 $, it follows
$ {\displaystyle \int_{\hat \Sigma} f(\mathbf x) \; d\hat\mu(\mathbf y, \mathbf x) < \frac{\epsilon}{2}} $.

We get then a contradiction. Therefore, $ \hat \mu(\hat R) = 0 $.
\end{proof}

We need also a result on numerical sequences.

\begin{3}\label{comportamentosequencias}
Consider a sequence $ \{a_j \} \subset \mathbb R $ for which is true
$$ \lim_{k \to \infty} \frac{1}{k} \sum_{j = 1}^{k} a_j = b. $$
Let $ R $ be a subset of the set of positive integers satisfying
$$ \lim_{k \to \infty} \frac{1}{k} \#\{j \in R: j \leq k \} > 0. $$
Then, for any $ \epsilon > 0 $ and any positive integer $ K $,
there exist $ k_1, k_2 \in R $ such that $ k_2 > k_1 \ge K $ and
$$ \left | \sum_{j = k_1 + 1}^{k_2} a_j - (k_2 - k_1)b \right | < \epsilon. $$
\end{3}

The previous lemma was used by Ma\~n\'e in \cite{Mane}. We can present now the following result.

\begin{4}
Suppose $ \sigma: \Sigma \to \Sigma $ is a transitive subshift of finite type. Let
$ A $ be a $\theta$-H\"older potential. Assume $ \hat \mu \in \text{\Large $\mathit m $}_A $ with $\hat \mu \circ \pi_1^{-1} $ ergodic.
Then $ \pi_1(\text{supp}(\hat \mu)) \subset \Omega(A) $.
\end{4}

\begin{proof}
It is enough to show that $ (\hat \mu \circ \pi_1^{-1})(\Omega(A)) = 1 $. Fix $ \epsilon > 0 $.
Denote by $ \Omega(A, \epsilon) $ the set of the points $ \mathbf x \in \Sigma $ for which we can find a path
$ \{(\mathbf y^0, \mathbf x^0), \ldots, (\mathbf y^{k - 1}, \mathbf x^{k - 1}) \} \in \mathcal P(\mathbf x, \mathbf x, \epsilon) $
satisfying
$$ \left | \sum_{j = 0}^{k - 1} (A - \beta_A)(\mathbf y^j, \mathbf x^j) \right | < \epsilon. $$
As $ \Omega(A) = \bigcap \Omega(A, 1/j) $, it is enough to show that $ (\hat \mu \circ \pi_1^{-1}) (\Omega(A, \epsilon)) = 1 $.

Suppose, however, that $ (\hat \mu \circ \pi_1^{-1}) (\pi_1(\text{supp}(\hat \mu)) - \Omega(A, \epsilon)) > 0 $.
Take an integer $ l > 0 $ sufficiently large in such way that $ 2 \lambda^l < \epsilon $. So
there exists $ \mathbf x \in \pi_1(\text{supp}(\hat \mu)) $ such that
$ (\hat \mu \circ \pi_1^{-1}) (D_l - \Omega(A, \epsilon)) > 0 $, where $ D_l $ is the open ball of radius
$ \lambda^l $ centered at the point $ \mathbf x $.

Thus, consider a point $\bar{\mathbf x} \in
\pi_1(\text{supp}(\hat \mu)) $ such that $$\lim_{k \to \infty}
\frac{1}{k} \#\{0 \leq j < k: \sigma^j(\bar{\mathbf x}) \in D_l -
\Omega(A, \epsilon) \} > 0. $$
Thanks to lemma~\ref{suporteholonomicas}, we can assume that, for every index $ j > 0 $, there
exists a point $ \bar{\mathbf y}^j \in \Sigma^* $ such that
$ (\bar{\mathbf y}^j,\sigma^j(\bar{\mathbf x})) \in \text{supp}(\hat \mu) $ and
$ \sigma^{j - 1}(\bar{\mathbf x}) = \tau_{\bar{\mathbf y}^j}(\sigma^j(\bar{\mathbf x})) $.

Being $ u \in C^0(\Sigma) $ an arbitrary sub-action for $ A $, from
proposition~\ref{contatomaximizantes} we get
$ A(\bar{\mathbf y}^j, \sigma^j(\bar{\mathbf x})) - \beta_A = u(\sigma^{j - 1}(\bar{\mathbf x})) - u(\sigma^j(\bar{\mathbf x})) $.
Define, finally,
$$ a_j = u(\sigma^{j - 1}(\bar{\mathbf x})) - u(\sigma^j(\bar{\mathbf x})) \; \text {and} \;
R = \{j: \sigma^j(\bar{\mathbf x}) \in D_l - \Omega(A, \epsilon) \}. $$

Using lemma~\ref{comportamentosequencias}, we obtain integers $ k_1, k_2 \in R $, with $ 1 \le k_1 < k_2 $, accomplishing
$$ \left | \sum_{j = k_1 + 1}^{k_2} (A - \beta_A)(\bar{\mathbf y}^j, \sigma^j(\bar{\mathbf x})) \right |
= \left | \sum_{j = k_1 + 1}^{k_2} a_j \right | < \epsilon. $$

However, once $ \sigma^{k_1}(\bar{\mathbf x}), \sigma^{k_2}(\bar{\mathbf x}) \in D_l $,
it follows that $ d(\sigma^{k_1}(\bar{\mathbf x}), \sigma^{k_2}(\bar{\mathbf x})) \le 2 \lambda^l $.
Therefore,
$ \{(\bar{\mathbf y}^{k_2}, \sigma^{k_2}(\bar{\mathbf x})), \ldots, (\bar{\mathbf y}^{k_1 + 1}, \sigma^{k_1 + 1}(\bar{\mathbf x})) \} \in
\mathcal P(\sigma^{k_2}(\bar{\mathbf x}), \sigma^{k_2}(\bar{\mathbf x}), \epsilon) $ yields
$ \sigma^{k_2}(\bar{\mathbf x}) \in \Omega(A, \epsilon) $. This is a contradiction because $ k_2 \in R $.

Hence, $ (\hat \mu \circ \pi_1^{-1}) (\Omega(A, \epsilon)) = 1 $.
\end{proof}

Remember that the addition of a constant does not change the role played by a sub-action.
Thus, the next proposition indicates a kind of rigidity created by the previous ergodic assumption.

\begin{4}
Consider a probability $ \hat \mu \in \text{\Large $\mathit m $}_A $ such that
$ \hat \mu \circ \pi_1^{-1} $ is ergodic. If $ u, u' \in C^0(\Sigma) $ are
sub-actions for $ A \in C^0(\hat \Sigma) $, then $ u - u' $ is
identically constant on $ \pi_1(\text{supp}(\hat \mu)) $.
\end{4}

\begin{proof}
Suppose $ \mathbf x \in \pi_1(\text{supp}(\hat \mu)) $. We can use
lemma~\ref{suporteholonomicas} in order to get a point $ (\bar{\mathbf y}, \bar{\mathbf x}) \in \text{supp}(\hat \mu) $
such that $ \mathbf x = \tau_{\bar{\mathbf y}} (\bar{\mathbf x}) $.
From proposition~\ref{contatomaximizantes}, we verify
$$ u(\bar{\mathbf x}) - u(\mathbf x) = \beta_A - A(\bar{\mathbf y}, \bar{\mathbf x}) = u'(\bar{\mathbf x}) - u'(\mathbf x). $$
So $ (u - u') (\mathbf x) = (u - u') (\bar{\mathbf x}) = (u - u') \circ \sigma (\mathbf x) $.
Therefore, we have $ u - u' = (u - u') \circ \sigma $ on $ \pi_1(\text{supp}(\hat \mu)) $.
As the probability $ \hat \mu \circ \pi_1^{-1} $ is ergodic,
it follows immediately that $ u - u' $ is constant on $\pi_1(\text{supp}(\hat \mu)) $.
\end{proof}

Let us consider again the transitivity hypothesis and assume $ A $ is
H\"older. Given $ u $ a sub-action for $ A $, let $ \mathbb M_A(u) $ be its corresponding contact locus.
Then, we claim that $ \Omega(A) \subset \pi_1(\mathbb M_A(u)) $. This is completely obvious when $ u $ is
a calibrated sub-action, because in such case $ \pi_1(\mathbb M_A(u)) = \Sigma $. Besides,
corollary~\ref{comportamentocalibradageral} tells us that every sub-action $ u \in C^0(\Sigma) $ for the potential $ A $
behaves as a calibrated sub-action on $ \Omega(A) $.

Therefore, the following inclusions are true
$$ \bigcup_{\stackrel{\hat \mu \in \text{\large $\mathit m $}_A}{\hat \mu \circ \pi_1^{-1} \text{ ergodic}}} \pi_1(\text{supp}(\hat \mu))
\; \subset \; \Omega(A) \; \subset
\bigcap_{\stackrel{u \in C^0(\Sigma)} {u \text{ sub-action}}} \pi_1(\mathbb M_A(u)). $$

In some situations for the standard model $ (X, T, \mathcal M_T) $, it
is known that, given a H\"older potential $ A $, a probability is $A$-maximizing
if, and only if, its support is contained in the set of non-wandering points
(with respect to $ A $). See, for instance, the case of expanding maps of the circle
in proposition 15.ii of \cite{CLT} and also the case of Anosov diffeomorphisms in lemmas 12 and 13 of
\cite{LT1}.

Hence, it is natural to ask: in order to verify that $ \hat \mu \in \text{\Large $\mathit m $}_A $,
it would be enough to check that $ \hat \mu \circ \pi_1^{-1} $ is ergodic and
$ \pi_1(\text{supp}(\hat \mu)) \subset \Omega(A) $? The answer is no.

Indeed, here is a counter-example. Take a potential $ A: \{0, 1 \}^{\mathbb Z} \to \mathbb R $ depending just on
three coordinates in such way that $ A(1,1|1) > A(s,s'|s '') $ whenever $ s + s' + s '' \le 2 $.
If we denote by $ \underline{ss'} $ either the periodic point $ (s, s', \ldots, s, s', \ldots) \in \Sigma $,
or the periodic point $ (\ldots, s, s', \ldots, s, s') \in \Sigma ^ * $, then we have
$ \delta_{(\underline{11},\underline{11})}, \delta_{(\underline{01},\underline{11})} \in \mathcal M_0 $
with $ \delta_{(\underline{11},\underline{11})} \circ \pi_1^{-1} = \delta_{\underline{11}} =
\delta_{(\underline{01},\underline{11})} \circ \pi_1^{-1} $. Nevertheless, observe that
$ \delta_{(\underline{11},\underline{11})} $ is a maximizing probability, but clearly
$ \delta_{(\underline{01},\underline{11})} \notin \text{\Large $\mathit m $}_A $.

The second inclusion above also bring us an interesting question:
what can be said about $ \pi_1(\mathbb M_A(u)) - \Omega(A) $? The
next proposition gives a partial answer.

\begin{4}\label{refinocontato}
Let $\sigma: \Sigma \to \Sigma $ be a transitive subshift of finite type and assume $ A \in
C^\theta(\hat \Sigma) $ is not cohomologous to a constant. Take
$ u \in C^0(\Sigma) $ an arbitrary sub-action for $ A $. Then, for
each positive integer $ k $, there exists a sub-action $ U_k \in C^0(\Sigma) $ satisfying
$$ \pi_1(\mathbb M_A(U_k)) \subset \bigcap_{j = 0}^{k - 1} \sigma^{-j}(\pi_1(\mathbb M_A(u))). $$
Moreover, if $ u $ is $\theta$-H\"older, then we can also take $ U_k $ as a $\theta$-H\"older function.
\end{4}

\begin{proof}
We begin with $ A^u = A + u \circ \pi_1 - u \circ \pi_1 \circ \hat \sigma^{-1} \leq \beta_A $.

Given $ k > 0 $ and $ \mathbf x \in \Sigma $, we call a path of size
$ k $ ending at the point $ \mathbf x $ any ordered sequence of
points $ (\mathbf y^0, \mathbf x^0), \ldots, (\mathbf y^{k - 1}, \mathbf x^{k - 1}) \in \hat \Sigma $
which verifies $ \mathbf x^0 = \mathbf x $ and $\mathbf x^{j + 1} = \tau_{\mathbf y^j}(\mathbf x^j) $ for $ 0 \leq j < k - 1 $.
Denote by  $ \mathcal P_k(\mathbf x) $ the set of such paths. Note that
$$ \sum_{j = 0}^{k - 1} A^{u}(\mathbf y^j, \mathbf x^j) \le k \beta_A $$
for $ \{(\mathbf y^0, \mathbf x^0), \ldots, (\mathbf y^{k - 1}, \mathbf x^{k - 1}) \} \in \mathcal P_k(\mathbf x) $.

Taking $ \{(\mathbf y^0, \sigma^{k-1}(\mathbf x)), (\mathbf y^1, \sigma^{k-2}(\mathbf x)), \ldots, (\mathbf y^{k-1},\mathbf x) \}
\in \mathcal P_k(\sigma^{k-1}(\mathbf x)) $, we have the identity
\begin{flushleft}
$ {\displaystyle \sum_{j = 0}^{k - 1} A(\mathbf y^j, \sigma^{k - 1 - j}(\mathbf x)) =} $ \\
\raggedleft {$ {\displaystyle = k A(\mathbf y^{k - 1}, \mathbf x) + \sum_{j = 0}^{k - 1} jA(\mathbf y^{j - 1}, \sigma^{k - j}(\mathbf x))
- \sum_{j = 0}^{k - 1} j A(\mathbf y^j, \sigma^{k - 1 - j}(\mathbf x))} $.}
\end{flushleft}

Now we define $ W: \Sigma \to \mathbb R$ in the following way
$$ W(\mathbf x) = \max_ {\{(\mathbf y^0, \sigma^{k-1}(\mathbf x)), \ldots, (\mathbf y^{k-1}, \mathbf x) \} \in \mathcal P_k(\sigma^{k-1}(\mathbf x))}
\left [\frac{1}{k}\sum_{j = 1}^{k - 1} jA(\mathbf y^{j - 1}, \sigma^{k - j}(\mathbf x)) \right]. $$
Once the correspondence $ \mathbf x \mapsto {\displaystyle \max_{y_0 = x_0} A(\mathbf y, \sigma(\mathbf x))} $ is
$\theta$-H\"older, the same is true for the function $ W $.

Fix a point $ (\mathbf y, \mathbf x) \in \hat \Sigma $. Then consider a path
$$ \{(\mathbf y^0, \sigma^{k-1}(\mathbf x)), \ldots, (\mathbf y^{k-2}, \sigma(\mathbf x)), (\mathbf y, \mathbf x) \} \in \mathcal P_k(\sigma^{k-1}(\mathbf x)) $$
accomplishing
$$ \frac{1}{k}\sum_{j = 1}^{k - 1} jA(\mathbf y^{j - 1}, \sigma^{k - j}(\mathbf x)) = W(\mathbf x). $$
Put $ \mathbf y^{k - 1} = \mathbf y $.
As $ \{(\mathbf y^1, \sigma^{k-2}(\mathbf x)), \ldots, (\mathbf y^{k-1}, \mathbf x) \} \in \mathcal P_{k-1}(\sigma^{k-1}(\tau_{\mathbf y}(\mathbf x))) $,
without difficulty we get
\begin{flushleft}
${\displaystyle A(\mathbf y, \mathbf x) + W(\mathbf x) - W(\tau_{\mathbf y}(\mathbf x))
\le} $\\
\centering {${\displaystyle \le A(\mathbf y^{k - 1}, \mathbf x) +
\frac{1}{k}\sum_{j = 0}^{k - 1} j A(\mathbf y^{j - 1}, \sigma^{k -
j}(\mathbf x))
- \frac{1}{k}\sum_{j = 0}^{k - 1} j A(\mathbf y^j, \sigma^{k - 1 - j}(\mathbf x)) =} $} \\
\raggedleft {${\displaystyle = \frac{1}{k}\sum_{j = 0}^{k - 1} A(\mathbf y^j, \sigma^{k - 1 - j}(\mathbf x))} $.}
\end{flushleft}
Therefore, if we denote $ U_k = W + k^{-1} S_k u $, we obtain
\begin{flushleft}
${\displaystyle A(\mathbf y, \mathbf x) + U_k(\mathbf x) -
U_k(\tau_{\mathbf y}(\mathbf x)) \le} $\\
\centering {${\displaystyle \le \frac{1}{k}\sum_{j = 0}^{k - 1}
A(\mathbf y^j, \sigma^{k - 1 - j}(\mathbf x)) +
\frac{1}{k} S_k u(\mathbf x) - \frac{1}{k}S_k u(\tau_{\mathbf y}(\mathbf x)) =} $} \\
\raggedleft {${\displaystyle = \frac{1}{k}\sum_{j = 0}^{k - 1}
A^u(\mathbf y^j, \sigma^{k - 1 - j}(\mathbf x)) \leq \beta_A} $.}
\end{flushleft}
Hence, $ U_k $ is a sub-action for the potential $ A $.

Let us check that such sub-action $ U_k $ accomplishes the claim of
the proposition. We just follow the itinerary of the construction of
$ U_k $ in the opposite direction. If $ \mathbf x \in \pi_1(\mathbb M_A(U_k)) $,
then there exists a path
$$ \{(\mathbf y^0, \sigma^{k-1}(\mathbf x)), \ldots, (\mathbf y^{k-1},\mathbf x) \} \in \mathcal P_k(\sigma^{k-1}(\mathbf x)) $$
such that
$$ \frac{1}{k}\sum_{j = 0}^{k - 1} A^u(\mathbf y^j, \sigma^{k - 1 - j}(\mathbf x)) = \beta_A, $$
which yields $ A^u(\mathbf y^j, \sigma^{k - 1 - j}(\mathbf x)) = \beta_A $. Thus, clearly
$ \sigma^{k - 1 - j}(\mathbf x) \in \pi_1(\mathbb M_A(u)) $ for all  $j \in \{0, \ldots, k - 1 \} $.
\end{proof}

The proof described above found inspiration in the strategy used
by Bousch in \cite{Bousch3}.

The previous proposition brings our attention to the following question:
does exist a non-calibrated sub-action? The answer is yes.

Under the same hypotheses of proposition~\ref{refinocontato}, assume that $ u \in C^\theta(\Sigma) $ is a calibrated sub-action.
Suppose yet the existence of a point $ (\mathbf y^0, \mathbf x^0) \in \hat \Sigma $
satisfying both $ {\displaystyle A(\mathbf y^0, \mathbf x^0) = \max_{y_0 = y_0^0} A(\mathbf y, \mathbf x^0)} $ and
$$ A(\mathbf y^0, \mathbf x^0) + u(\mathbf x^0) - u(\tau_{\mathbf y^0}(\mathbf x^0)) < \beta_A. $$
(These assumptions are obviously verified by any potential $ A \in C^\theta(\Sigma) $ not cohomologous to a constant.)
We claim that the function $ U \in C^\theta(\Sigma) $ defined by
$$ U(\mathbf x) = \frac{1}{2}[u(\sigma(\mathbf x)) + u(\mathbf x)] + \frac{1}{2} \max_{y_0 = x_0} A(\mathbf y, \sigma(\mathbf x)) $$
is a sub-action for $ A $ which is not calibrated. Indeed, the function $ U $
is nothing else that the sub-action $ U_2 $ described in the proof of the previous proposition.
Moreover, note that, for all $ \mathbf y \in \Sigma_{\tau_{\mathbf y^0}(\mathbf x^0)}^* $,
\begin{flushleft}
$A(\mathbf y,\tau_{\mathbf y^0}(\mathbf x^0)) + U(\tau_{\mathbf y^0}(\mathbf x^0)) -
U(\tau_{\mathbf y}(\tau_{\mathbf y^0}(\mathbf x^0))) \le $\\
\centering {$ {\displaystyle \le \frac{1}{2}
[ A(\mathbf y, \tau_{\mathbf y^0}(\mathbf x^0)) + u(\tau_{\mathbf y^0}(\mathbf x^0)) - u(\tau_{\mathbf y}(\tau_{\mathbf y^0}(\mathbf x^0)))] \; +} $} \\
\raggedleft {${\displaystyle + \; \frac{1}{2} [A(\mathbf y^0, \mathbf x^0) + u(\mathbf x^0) - u(\tau_{\mathbf y^0}(\mathbf x^0))]} <
\beta_A $,}
\end{flushleft}
therefore $ \tau_{\mathbf y^0}(\mathbf x^0) \notin \pi_1(\mathbb M_A(U)) $.

A deeper study of non-calibrated sub-actions is the aim of a subsequent paper \cite{GLT}.
Finally, we would like to mention that the possibility of adapting our holonomic setting 
to the case of iterated function systems has been recently announced \cite{Oliveira}.

\end{document}